Florentin Smarandache

# Neutrosophic Precalculus and Neutrosophic Calculus

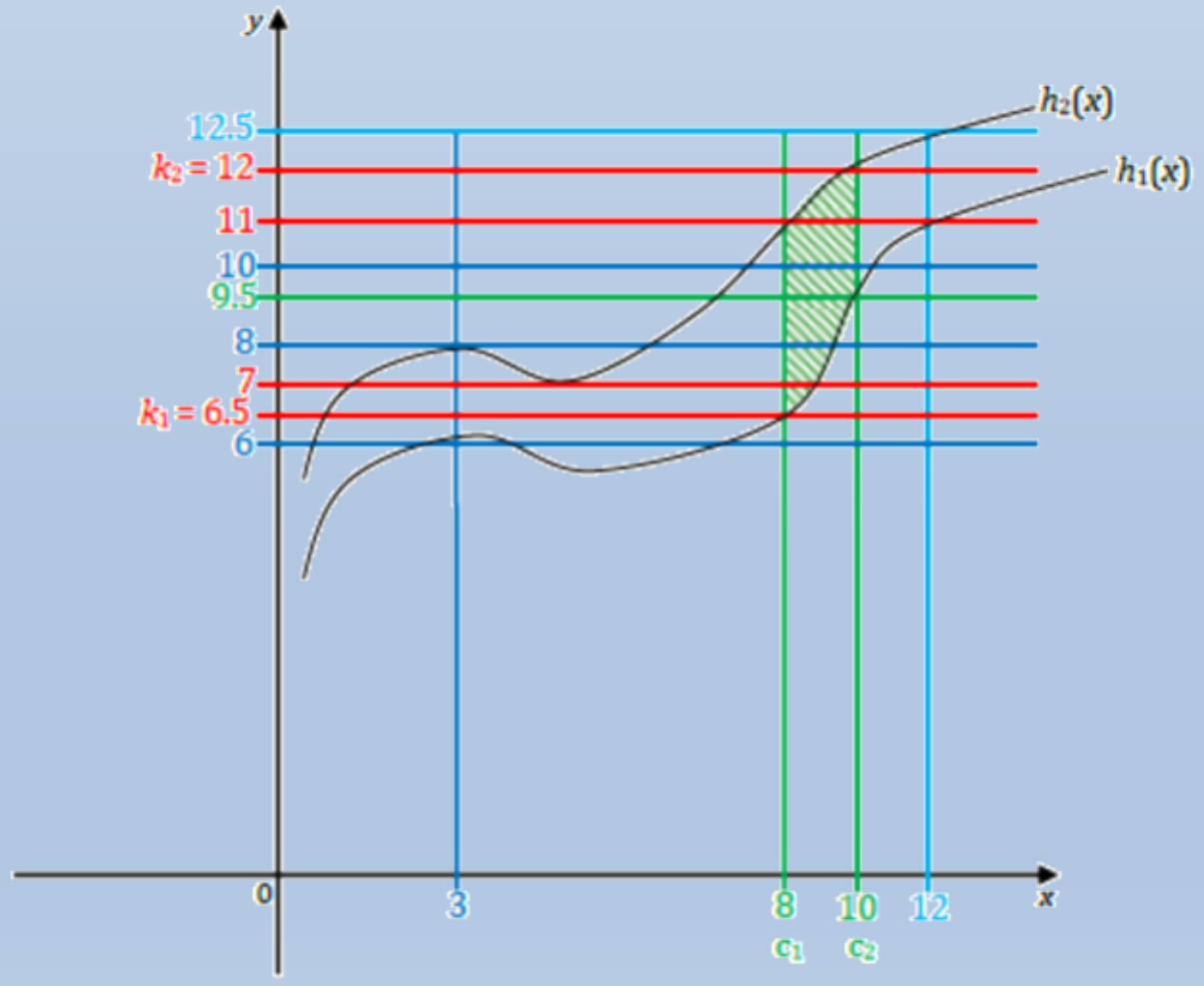

Florentin Smarandache

# Neutrosophic Precalculus and Neutrosophic Calculus

**2015**

Florentin Smarandache

# Neutrosophic Precalculus

and

# Neutrosophic Calculus

*On the frontcover:* Example for the Neutrosophic
Intermediate Value Theorem

# Contents

# I. Introductory Remarks

# I.1. Overview

*Neutrosophy* means the study of ideas and notions that are not true, nor false, but in between (i.e. neutral, indeterminate, unclear, vague, ambiguous, incomplete, contradictory, etc.).

Each field has a neutrosophic part, i.e. that part that has indeterminacy. Thus, there were born the neutrosophic logic, neutrosophic set, neutrosophic probability, neutrosophic statistics, neutrosophic measure, neutrosophic precalculus, neutrosophic calculus, etc.

There exist many types of indeterminacies – that's why neutrosophy can be developed in many different ways.

# I.2. Preliminary

The first part of this book focuses on *Neutrosophic Precalculus*, which studies the neutrosophic functions. A *Neutrosophic Function* $f: A \rightarrow B$ is a function which has some indeterminacy, with respect to its domain of definition, to its range, or to its relationship that associates elements in $A$ with elements in $B$.

As particular cases, we present the *neutrosophic exponential function* and *neutrosophic logarithmic function*. The *neutrosophic inverse function* is the inverse of a neutrosophic function.

A *Neutrosophic Model* is, in the same way, a model with some indeterminacy (vagueness, unsureness, ambiguity, incompleteness, contradiction, etc.).

*

The second part of the book focuses on *Neutrosophic Calculus*, which studies the neutrosophic limits, neutrosophic derivatives, and neutrosophic integrals.

*

We introduce for the first time the notions of *neutrosophic mereo-limit, mereo-continuity, mereo-derivative,* and *mereo-integral*, [1] besides the classical

[1] From the Greek μερος, 'part'. It is also used to define the theory of the relations of part to whole and the relations of part to part within a whole (mereology), started by Leśniewski, in "Foundations of the General Theory of Sets" (1916) and "Foundations of Mathematics" (1927–1931), continued by Leonard and Goodman's "The Calculus of Individuals" (1940),

definitions of limit, continuity, derivative, and integral respectively.

*

The *Neutrosophic Precalculus* and *Neutrosophic Calculus* can be developed in many ways, depending on the types of indeterminacy one has and on the method used to deal with such indeterminacy.

In this book, we present a few examples of indeterminacies and several methods to deal with these specific indeterminacies, but many other indeterminacies there exist in our everyday life, and they have to be studied and resolved using similar of different methods. Therefore, more research has to be done in the field of neutrosophics.

# I.3. Distinctions among Interval Analysis, Set Analysis, and Neutrosophic Analysis

## Notation

In this book we consider that an interval *[a, b] = [b, a]* in the case when we do not know which one between *a* and *b* is bigger, or for the case when the interval has varying left and right limits of the form *[f(x), g(x)]*, where for certain *x*'s one has *f(x) < g(x)* and for other *x*'s one has *f(x) > g(x).*

## Interval Analysis

In **Interval Analysis** (or Interval Arithmetic) one works with intervals instead of crisp numbers. Interval analysis is intended for rounding up and down errors of calculations. So an error is bounding by a closed interval.

## Set Analysis

If one replaces the closed intervals (from interval analysis) by a set, one get a **Set Analysis** (or Set Arithmetic).

For example, the set-argument set-value function:

$h: P(\mathrm{R}) \rightarrow P(\mathrm{R})$, (1)

where *P*(R) is the power set of R (the set of all real numbers),

$h(\{1, 2, 3\}) = \{7, 9\}$, $h([0, 1]) = (6, 8)$, $h(-3) = \{-1, -2\} \cup (2.5, 8]$, $h([x, x^2] \cup [-x^2, x]) = 0$. (2)

Set analysis is a generalization of the interval analysis.

## Distinctions among Interval Analysis, Set Analysis, and Neutrosophic Analysis

**Neutrosophic Analysis** (or Neutrosophic Arithmetic) is a generalization of both the interval analysis and set analysis, because neutrosophic analysis deals with all kind of sets (not only with intervals), and also considers the case when there is some indeterminacy (with respect to the sets, or with respect to the functions or other notions defined on those sets).

If one uses sets and there is no indeterminacy, then neutrosophic analysis coincides with the set analysis.

If instead of sets, one uses only intervals and there is no indeterminacy, then neutrosophic analysis coincides with interval analysis.

If there is some indeterminacy, no matter if using only intervals, or using sets, one has neutrosophic analysis.

## Examples of Neutrosophic Analysis

Neutrosophic precalculus and neutrosophic calculus are also different from set analysis, since they use indeterminacy.

As examples, let's consider the *neutrosophic functions*:

$f_1(0 \text{ or } 1) = 7$ (indeterminacy with respect to the argument of the function),

i.e. we are not sure if $f_1(0) = 7$ or $f_1(1) = 7$. (3)

Or

$f_2(2) = 5 \text{ or } 6$ (indeterminacy with respect to the value of the function),

so we are not sure if $f_2(2) = 5$ or $f_2(2) = 6$. (4)

Or even more complex:

$f_3(-2 \text{ or } -1) = -5 \text{ or } 9$ (indeterminacy with respect with both the argument and the value of the function),

i.e. $f_3(-2) = -5$, or $f_3(-2) = 9$, or $f_3(-1) = -5$, or $f_3(-1) = 9$. (5)

And in general:

$f_{m,n}(a_1 \text{ or } a_2 \text{ or } \dots \text{ or } a_m) = b_1 \text{ or } b_2 \text{ or } \dots \text{ or } b_n$. (6)

These functions, containing such indeterminacies, are different from set-valued vector-functions.

## Examples in Set Analysis

For example $f_1$: R → R is different from the set-argument function:

$g_1: R^2 \rightarrow R$, where $g_1(\{0, 1\}) = 7$. (7)

Also, $f_2$: R → R is different from the set-value function

$g_2: R \rightarrow R^2$, where $g_2(2) = \{5, 6\}$. (8)

Similarly, $f_3$: R → R is different from the set-argument set-value function

$g_3: R^2 \rightarrow R^2$, where $g_3(\{-2, -1\}) = \{-5, 9\}$. (9)

And in the general case, $f_{m,n}$: R → R is different from the set-argument set-value function

$g_{m,n} : R^m \rightarrow R^n$,

where $g_{m,n}(\{a_1, a_2, \dots, a_m\}) = \{b_1, b_2, \dots, b_n\}$. (10)

It is true that any set can be enclosed into a closed interval, yet by working with larger intervals than narrow sets, the result is rougher, coarser, and more inaccurate.

Neutrosophic approach, by using smaller sets included into intervals, is more refined than interval analysis.

Neutrosophic approach also uses, as particular cases, open intervals, and half-open half-closed intervals.

## Examples in Interval Analysis

Also, neutrosophic analysis deals with sets that have some indeterminacy: for example we know that an element *x(t,i,f)* only partially belongs to a set *S*, and partially it does not belong to the set, while another part regarding the appurtenance to the set is indeterminate.

Or we have no idea if an element *y(0,1,0)* belongs or not to the set (complete indeterminacy).

Or there is an element that belongs to the set, but we do not know it.

Interval analysis and set analysis do not handle these.

Let's consider an interval L = $[0, 5_{(0.6,\ 0.1,\ 0.3)}[$, where the number 5(0.6, 0.1, 0.3) only partially (0.6) belongs to the interval L, partially doesn't belong (0.3), and its appurtenance is indeterminate (0.1). We should observe that L ≠ [0, 5] and L ≠ [0, 5). Actually, L is in between them:

$$[0, 5) \subset L \subset [0, 5], \tag{11}$$

since the element 5 does not belong to [0, 5), partially belong to $[0, 5_{(0.6,\ 0.1,\ 0.3)}[$, and certainly belongs to [0, 5]. So, the interval L is part of neutrosophic analysis, not of interval analysis.

Now, if one considers the functions:

$$k_1( [0, 5] ) = [-4, 6], \text{ or } k_2( [-2, -4] ) = [0, 5], \tag{12}$$

then $k_1$ and $k_2$ belong to the interval analysis.

But if we take

$k_3([0, 5_{(0.6,\ 0.1,\ 0.3)}[)=[-4, 6]$, or $k_4([-2, -4])=[0, 5_{(0.6,0.1,0.3)}[$,

then $k_3$ and $k_4$ belong to neutrosophic analysis. (13)

A *Neutrosophic Function* $f: A \to B$ is a function, which has some indeterminacy, with respect to its domain of definition, to its range, to its relationship that associates elements in $A$ with elements in $B$ -- or to two or three of the above situations.

Interval Analysis studies only functions defined on intervals, whose values are also intervals, but have no indeterminacy.

Therefore, neutrosophic analysis is more general than interval analysis. Also, neutrosophic analysis deals with indeterminacy with respect to a function argument, a function value, or both.

For example, the neutrosophic functions:

$$e: \mathbb{R} \cup \{I\} \to \mathbb{R} \cup \{I\}, e(2+3I) = 7-6I \tag{14}$$

where $I$ = indeterminacy.

$$f: \mathbb{R} \to \mathbb{R}, f(4 \text{ or } 5) = 7; \tag{15}$$

$$g: \mathbb{R} \to \mathbb{R}, g(0) = -2 \text{ or } 3 \text{ or } 7; \tag{16}$$

$$h: \mathbb{R} \to \mathbb{R}, h(-1 \text{ or } 1) = 4 \text{ or } 6 \text{ or } 8; \tag{17}$$

$k: \mathbb{R} \to \mathbb{R}, k(x) = x \text{ and } -x$ (which fails the classical vertical line test for a curve to be a classical function); thus $k(x)$ is not a function from a classical point of view, but it is a neutrosophic function); (18)

$$l: \mathbb{R} \to \mathbb{R}, l(-3) = \text{maybe } 9. \tag{19}$$

One has:

Interval Analysis ⊂ Set Analysis ⊂ Neutrosophic Analysis.

## Inclusion Isotonicity

Inclusion isotonicity of interval arithmetic also applies to set analysis and neutrosophic analysis. Hence, if $\odot$ stands for set addition, set subtraction, set multiplication, or set division, and A, B, C, D are four sets such that: $A \subseteq C$ and $B \subseteq D$, then

$$A \odot B \subseteq C \odot D. \tag{20}$$

The proof is elementary for set analysis:

Let $x \in A \odot B$, then there exists $a \in A$ and $b \in B$ such that $x = a \odot b$.

But $a \in A$ and $A \subseteq C$ means that $a \in C$ as well.

And similarly, $b \in B$ and $B \subseteq D$ means that $b \in D$ as well.

Whence, $x = a \odot b \in C \odot D$ too.

The proof for neutrosophic analysis is similar, but one has to consider one of the neutrosophic inclusion operators; for example as follows for crisp neutrosophic components *t, i, f*:

a neutrosophic set *M* is included into a neutrosophic set *N* if,

for any element $x(t_M,i_M,f_M) \in M$ one has $x(t_n,i_n,f_n) \in N$, with $t_M \leq t_N$, $i_M \geq i_N$, and $f_M \geq f_N$.

## Conclusion

This research is in the similar style as those on neutrosophic probability (2013) and neutrosophic statistics (2014) from below.

# I.4. Indeterminate Elementary Geometrical Measurements

The mathematics of indeterminate change is the *Neutrosophic Calculus*.

Indeterminacy means imprecise, unclear, vague, incomplete, inconsistent, contradictory information. While classical calculus characterizes the dynamicity of our world, neutrosophic calculus characterizes the indeterminate (neutrosophic) dynamicity. Classical calculus deals with notions (such as slope, tangent line, arc length, centroid, curvature, area, volume, as well as velocity, and acceleration) as exact measurements, but in many real-life situations one deals with approximate measurements.

*Neutrosophic Precalculus* is more static and is referred to ambiguous staticity.

In neutrosophic calculus, we deal with notions that have some indeterminacy. Moreover, indeterminacy, unfortunately, propagates from one operation to the other.

In an abstract idealist world, there are perfect objects and perfect notions that the classical calculus uses.

For example, the curvature of perfect circle of radius $r > 0$ is a constant number [equals to $1/r$], but for an imperfect circle its curvature may be an interval [included in $(1/r-\varepsilon, 1/r+\varepsilon)$, which is a neighborhood of the number 1/r, with $\varepsilon > 0$ a tiny number].

A perfect right triangle with legs of *1 cm* and *2 cm* has its hypotenuse equals to $\sqrt{5}$ cm. However, in our imperfect world, we cannot draw a segment of line whose length be equal of exactly $\sqrt{5}$ $cm$, since $\sqrt{5}$ is an irrational number that has infinitely many decimals, we need to approximate it to a few decimals: $\sqrt{5} = 2.23606797 \ldots$

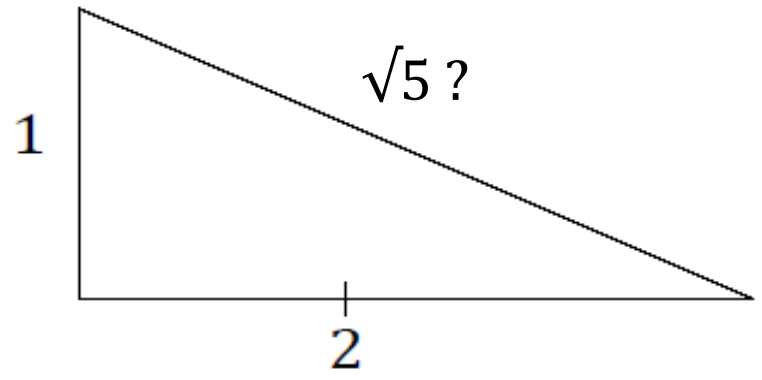

*Figure 1.*

The area of a perfect ellipses is $A = \pi ab$, where $2a$ and $2b$, with $a > b$, are its major and minor axes respectively. However, we cannot represent it exactly since $\pi$ is a transcendental number (i.e. it is not a solution of any polynomial equations with rational coefficients), and it has infinitely many decimals. If $a = 2\ cm$ and $b = 1\ cm$, then the area of the ellipse is $A = 2\pi = 6.2831 \ldots cm^2$.

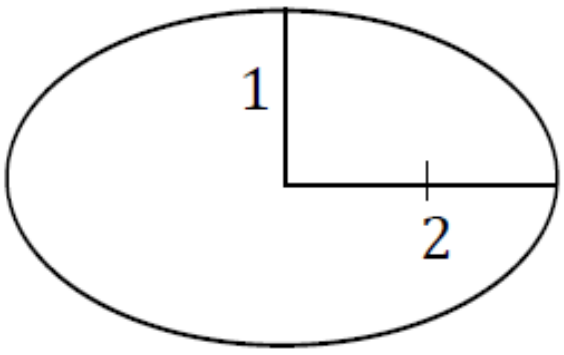

*Figure 2.*

but we can exactly comprise this area inside of this ellipse, since 6.2831 ... is not an exact number. We only work with approximations (imprecisions, indeterminations).

Similarly, for the volume of a perfect sphere $V = \frac{4}{3}\pi r^3$ where its radius is $r$. If $r = 1$ cm, then $V = \frac{4}{3}\pi = 4.1887 \ldots \text{cm}^3$ which is a transcendental number and has infinitely many decimals. Thus, we are not able to exactly have the volume of the below sphere,

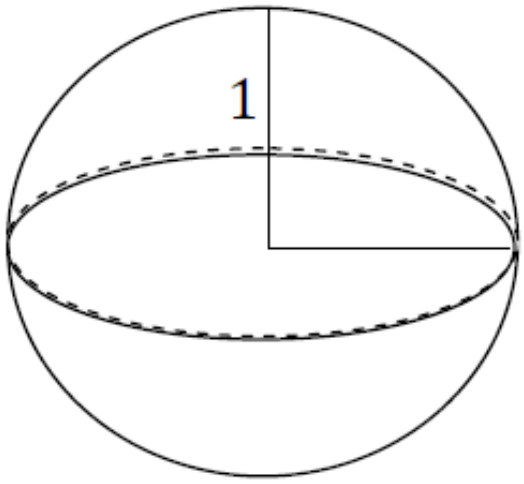

*Figure 3.*

equals to $4.1887 \ldots cm^3$.

# I.5. Indeterminate Physical Laws

Neutrosophy has also applications in physics, since many physical laws are defined in strictly closed systems, i.e. in idealist (perfect) systems[2], but such "perfect" system do not exist in our world, we deal only with approximately closed system, which makes room for using the neutrosophic (indeterminate) theory. Therefore, a system can be *t%* closed (in most cases $t < 100$), *i%* indeterminate with respect to closeness or openness, and *f%* open.

Therefore, a theoretical physical law ($L$) may be true in our practical world in less than *100%*, hence the law may have a small percentage of falsehood, and another small percentage of indeterminacy (as in neutrosophic logic).

Between the validity and invalidity of a theoretical law (idea) in practice, there could be included multiple-middles, i.e. cases where the theoretical law (idea) is partially valid and partially invalid.

[2] Fu Yuhua, "Pauli Exclusion Principle and the Law of Included Multiple-Middle", in Neutrosophic Sets and Systems, Vol. 6, 2014.

# II. Neutrosophic Precalculus

## II.1. Algebraic Operations with Sets

Let $S$ and $T$ be two sets, and $\alpha \in \mathbb{R}$ a scalar. Then:

$$\alpha \cdot S = \{\alpha \cdot s | s \in S\}; \quad (21)$$

$$S + T = \{s + t | s \in S, t \in T\}; \quad (22)$$

$$S - T = \{s - t | s \in S, t \in T\}; \quad (23)$$

$$S \cdot T = \{s \cdot t | s \in S, t \in T\}; \quad (24)$$

$$\frac{S}{T} = \left\{\frac{s}{t} | s \in S, t \in T, t \neq 0\right\}. \quad (25)$$

## II.2. Neutrosophic Subset Relation

A *Neutrosophic Subset Relation* $r$, between two sets $A$ and $B$, is a set of ordered pairs of the form $(S_A, S_B)$, where $S_A$ is a subset of $A$, and $S_B$ a subset of $B$, with some indeterminacy.

A neutrosophic relation $r$, besides *sure ordered pairs* $(S_A, S_B)$ that 100% belong to $r$, may also contains *potential ordered pairs* $(S_C, S_D)$, where $S_C$ is a subset of $A$, and $S_D$ a subset of $B$, that might be possible to belong to $r$, but we do not know in what degree, or that partially belong to $r$ with the neutrosophic value $(T, I, F)$, where $T < 1$ means degree of appurtenance to $r$, $I$ means degree of indeterminate appurtenance, and $F$ means degree of non-appurtenance.

*Example*:

$r: \{0, 2, 4, 6\} \rightarrow \{1, 3, 5\}$

$$r = \begin{Bmatrix} (\{0, 2\}, \{1, 3\}), (\{4, 6\}, \{5\}), \\ (\{6\}, \{1, 5\})_{(0.7, 0.1, 0.1)}, (\{2, 6\}, \{3, 5\})_? \end{Bmatrix} \quad (26)$$

where $(\{0, 2\}, \{1, 3\})$ and $(\{4, 6\}, \{5\})$ for sure belong to $r$; while $(\{6\}, \{1, 5\})$ partially belongs to $r$ in a percentage of 70%, 10% is its indeterminate appurtenance, and 10% doesn't belong to $r$;
and $(\{2, 6\}, \{3, 5\})$ is also potential ordered pairs (it might belong to $r$, but we don't know in what degree).

# II.3. Neutrosophic Subset Function

A *Neutrosophic Subset Function* $f: \mathcal{P}(A) \to \mathcal{P}(B)$, is a neutrosophic subset relation such that if there exists a subset $S \subseteq A$ with $f(s) = T$, and $f(s) = T_2$, then $T_1 \equiv T_2$. (This is the *(Neutrosophic) Vertical Line Test* extended from crisp to set-values.)

As a particular case, a *Neutrosophic Crisp Relation* between two sets $A$ and $B$ is a classical (crisp) relation that has some indeterminacy.

A neutrosophic crisp relation may contain, besides the *classical sure ordered pairs* $(a, b)$, with $a \in A$ and $b \in B$, also *potential ordered pairs* $(c, d)$, with $c \in A$ and $d \in B$ meaning that we are not sure if there is or there is not a relation between $c$ and $d$, or there is a relation between $c$ and $d$, but in a percentage strictly less then *100%*.

For example, the neutrosophic relation:

$$r: \{1, 2, 3, 4\} \to \{5, 6, 7, 8, 9\} \quad (27)$$

defined in set notation as:

$$\{(1, 5), (2, 6), (3, 7)_{[0.6, 0.1, 0.2]}, (3, 8)_?, (4, 9)_?\}$$

where the ordered pairs $(1, 5), (2, 6), (3, 7)$ for sure (*100%* belong to $r$), while $(3, 7)$ only *60%* belongs to $r$, *10%* the appurtenance is indeterminate, and *30%* it does not belong to $r$ [as in neutrosophic set], while about the ordered pairs $(3, 8)$ and $(4, 9)$ we do not know their appurtenance to $r$ (but it might be possible).

Another definition, in general, is:

A *Neutrosophic Relation* $r: A \to B$ is formed by any connections between subsets and indeterminacies in $A$ with subsets and indeterminacies in $B$.

It is a double generalization of the classical relation; firstly, because instead of connecting elements in $A$ with elements in $B$, one connects subsets in $A$ with subsets in $B$; and secondly, because it has some indeterminacies, or connects indeterminacies, or some connections are not well-known.

A neutrosophic relation, which is not a neutrosophic function, can be restrained to a neutrosophic function in several ways.

For example, if $r(S) = T_1$ and $r(S) = T_2$, where $T_1 \neq T_2$, we can combine these to:

- either $f(S) = T_1$ and $T_2$,
- or $f(S) = T_1$ or $T_2$,
- or $f(S) = \{T_1, T_2\}$,

which comply with the definition of a neutrosophic function.

# II.4. Neutrosophic Crisp Function

A *Neutrosophic Crisp Function* $f: A \rightarrow B$ is a neutrosophic crisp relation, such that if there exists an element $a \in A$ with $f(a) = b$ and $f(a) = c$, where $b, c \in B$, then $b \equiv c$. (This is the classical *Vertical Line Test*.)

## II.5. General Neutrosophic Function

A *General Neutrosophic Function* is a neutrosophic relation where the vertical line test (or the vertical subset-line text) does not work. But, in this case, the general neutrosophic function coincides with the neutrosophic relation.

# II.6. Neutrosophic (Subset or Crisp) Function

A *neutrosophic (subset or crisp) function* in general is a function that has some indeterminacy.

## Examples

1. $$f: \{1, 2, 3\} \rightarrow \{4, 5, 6, 7\} \quad (28)$$
$f(1) = 4, f(2) = 5$, but $f(3) = 6$ or 7 [we are not sure].

If we consider a *neutrosophic diagram representation* of this neutrosophic function, we have:

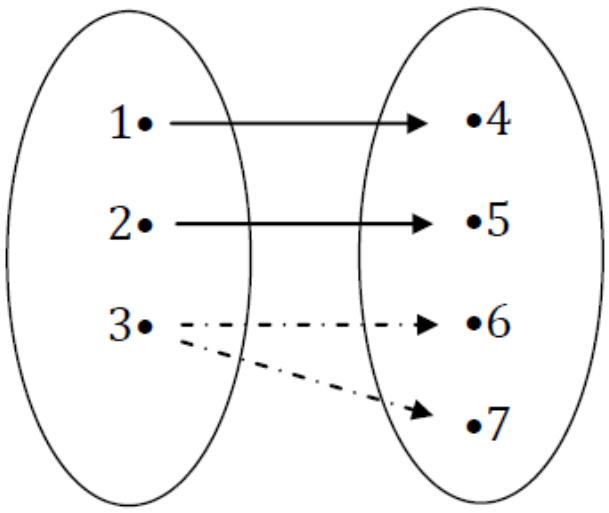

*Diagram 1.* Neutrosophic Diagram Representation.

The dotted arrows mean that we are not sure if the element 3 is connected to the element 6, or if 3 is connected to 7.

As we see, this neutrosophic function is not a function in the classical way, and it is not even a relationship in a classical way.

If we make a *set representation* of this neutrosophic function, we have:

$$\{(1,4), (2,5), \boxed{(3,6)}_?, \boxed{(3,7)}_?\}$$

where the dotted borders mean we are nou sure if they belong or not to this set. Or we can put the pairs (3, 6) and (3, 7) in red color (as warning).

In table representation, we have:

| x | 1 | 2 | 3 ? | 3 ? |
|---|---|---|---|---|
| f(x) | 4 | 5 | 6 | 7 |

*Table 1.*

where about the red color numbers we are not sure.

Similarly, for a *graph representation*:

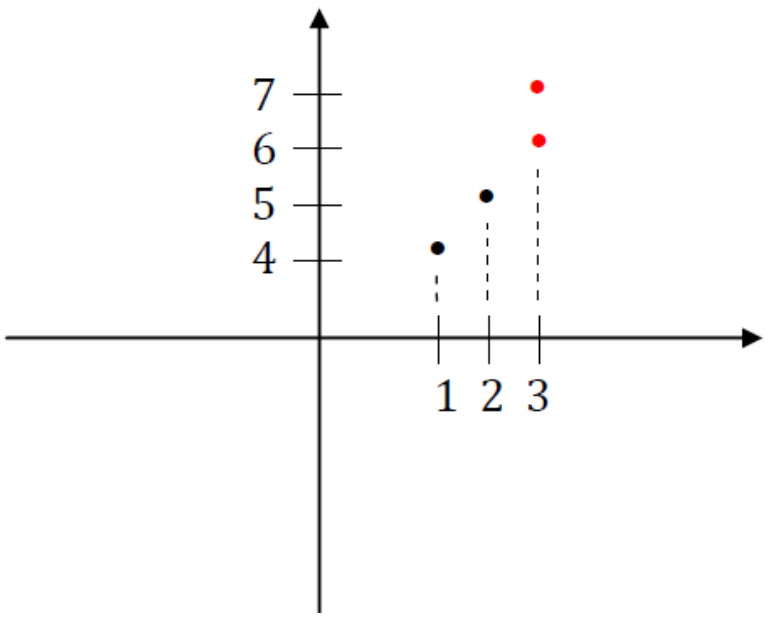

*Graph 1.*

Or, modifying a little this example, we might know, for example, that *3* is connected with *7* only partially, i.e. let's say $(3, 7)_{(0.6, 0.2, 0.5)}$ which means that 60% 3 is connected with 7, 20% it is not clear if connected or non-connected, and 50% 3 is not connected with 7.

The sum of components 0.6 + 0.2 + 0.5 = 1.3 is greater than 1 because the three sources providing information about connection, indeterminacy, non-connection respectively are independent, and use different criteria of evaluation.

2. We modify again this neutrosophic function as follows:

$$g: \{1,2,3\} \to \{4,5,6,7\}, \tag{29}$$

$g(1) = 4, g(2) = 5$, but $g(3) = 6$ and 7.

The neutrosophic function $g$ is not a function in the classical way (since it fails the vertical line test at $x = 3$), but it is a relationship in the classical way.

Its four representations are respectively:

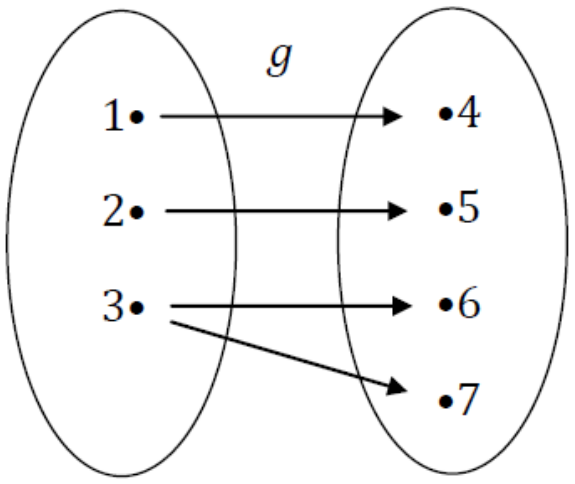

*Diagram 2.*

$$\{(1,4),(2,5),(3,6),(3,7)\}$$

| x | 1 | 2 | 3 | 3 |
|---|---|---|---|---|
| f(x) | 4 | 5 | 6 | 7 |

*Table 1.*

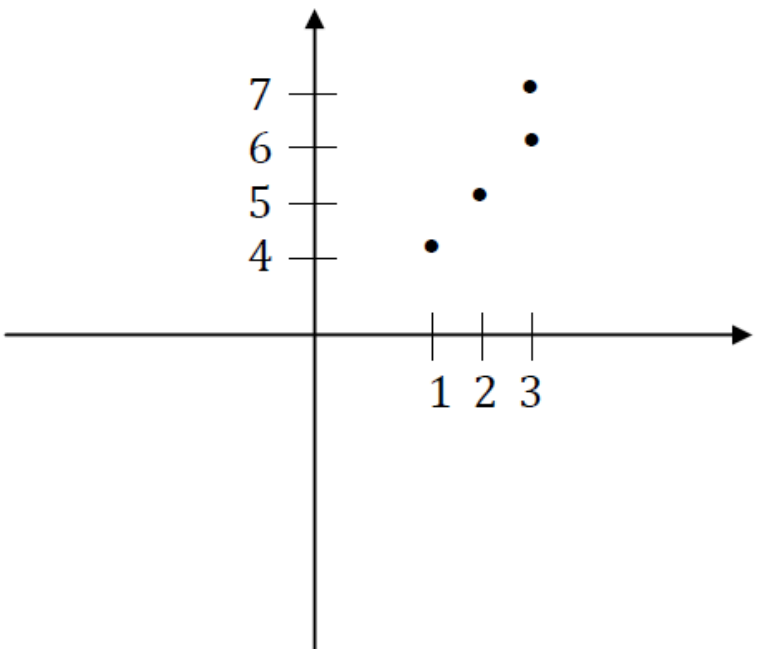

*Graph 2.*

Yet, if we redesign $g$ as

$G: \{1, 2, 3\} \to \mathcal{P}(\{4, 5, 6, 7\})$, (30)

$G(1) = 4$, $G(2) = 5$, and $G(3) = \{6, 7\}$,

then $G$ becomes a classical set-valued function.

3. Let's consider a different style of neutrosophic function:

$$h: \mathbb{R} \to \mathbb{R} \quad (31)$$

$$h(x) \in [2, 3], \text{ for any } x \in \mathbb{R}.$$

Therefore, we know about this function only the fact that it is bounded by the horizontal lines $y = 2$ and $y = 3$:

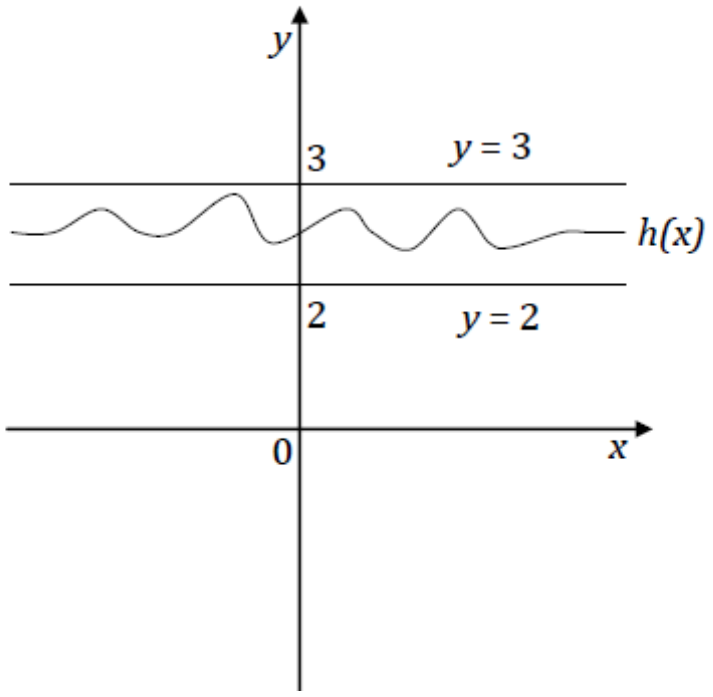

*Graph 3.*

4. Similarly, we modify $h(\bullet)$ and get a constant neutrosophic function (or thick function):

$l: \mathbb{R} \to \mathcal{P}(\mathbb{R})$ (32)

$l(x) = [2, 3]$ for any $x \in \mathbb{R}$,

where $\mathcal{P}(\mathbb{R})$ is the set of all subsets of $\mathbb{R}$.

For ex., $l(7)$ is the vertical segment of line [2, 3].

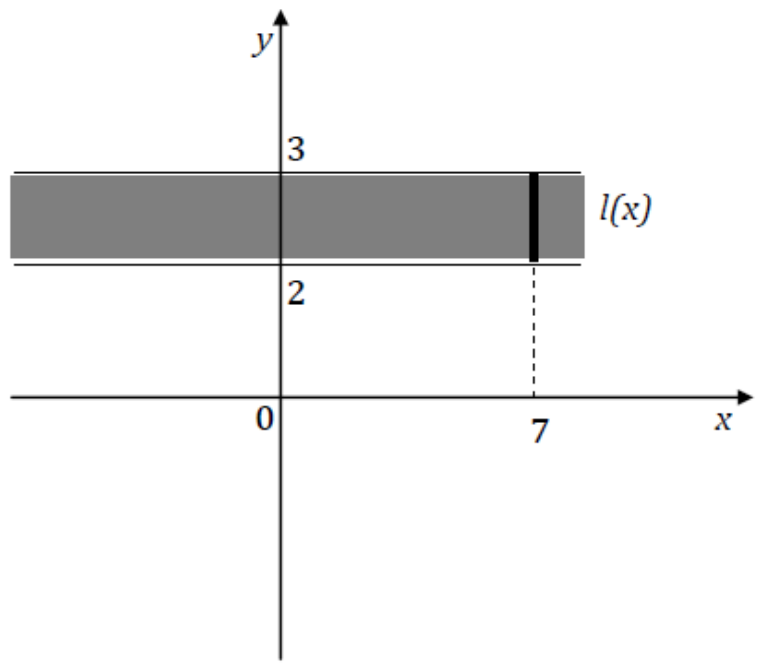

*Graph 4.*

5. A non-constant neutrosophic thick function:

$$k: \mathbb{R} \rightarrow \mathcal{P}(\mathbb{R}) \quad (33)$$

$$k(x) = [2x, 2x + 1]$$

whose graph is:

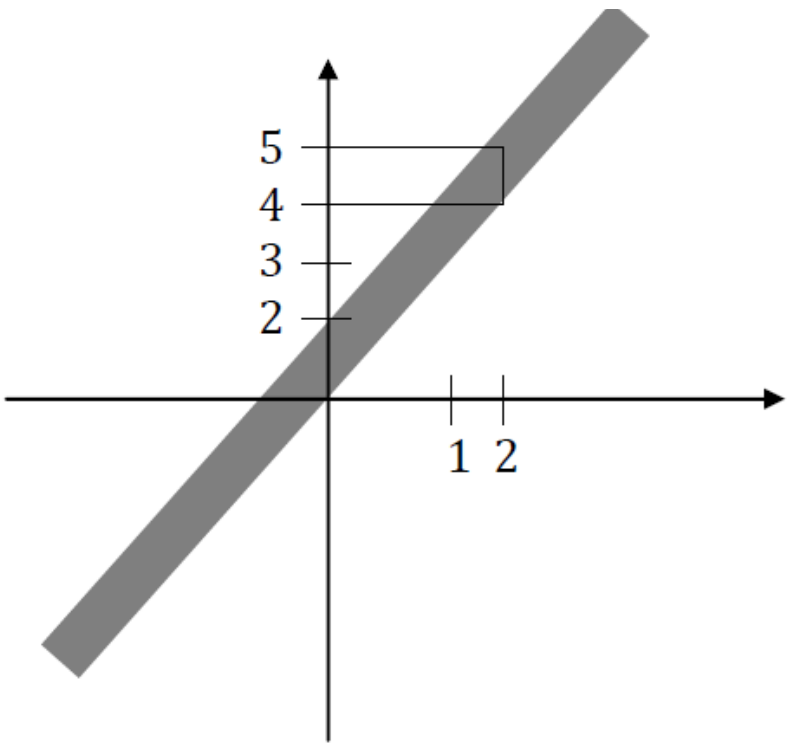

*Graph 5.*

For example:

$$k(2) = [2(2), 2(2) + 1] = [4, 5].$$

6. In general, we may define a neutrosophic thick function as:

$$m: \mathbb{R} \to \mathcal{P}(\mathbb{R}) \quad (34)$$

$$m(x) = [m_1(x_1)m_2(x)]$$

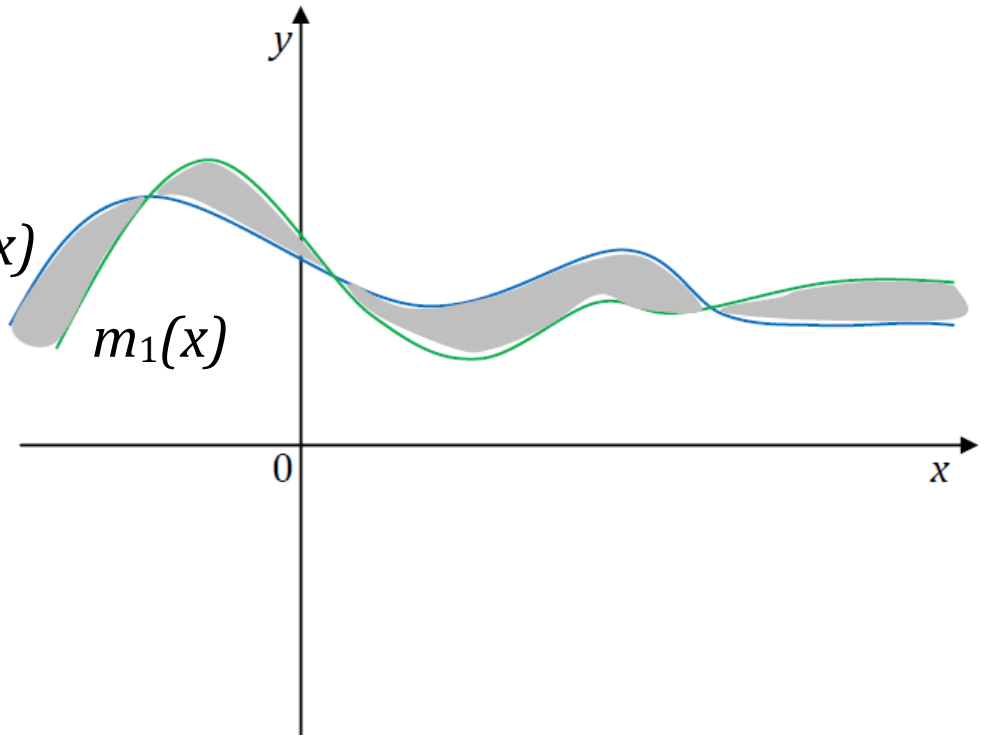

*Graph 6.*

and, of course, instead of brackets we may have an open interval $(m_1(x), m_2(x))$, or semi-open/semi-close intervals $(m_1(x), m_2(x)]$, or $[m_1(x),\ m_2(x)]$.

For example, $m(0) = [m_1(0), m_2(0)]$, the value of neutrosophic function $m(x)$ and a vertical segment of line.

These examples of thick (neutrosophic) functions are actually classical surfaces in $\mathbb{R}^2$.

7. Example of neutrosophic piece-wise function:

$$s: \mathbb{R} \to \mathcal{P}(\mathbb{R}) \quad (35)$$

$$s(x) = \begin{cases} [s_1(x), s_2(x)], \text{for } x \leq 3; \\ (s_3(x), s_4(x), \text{for } x > 3; \end{cases}$$

with the neutrosophic graph:

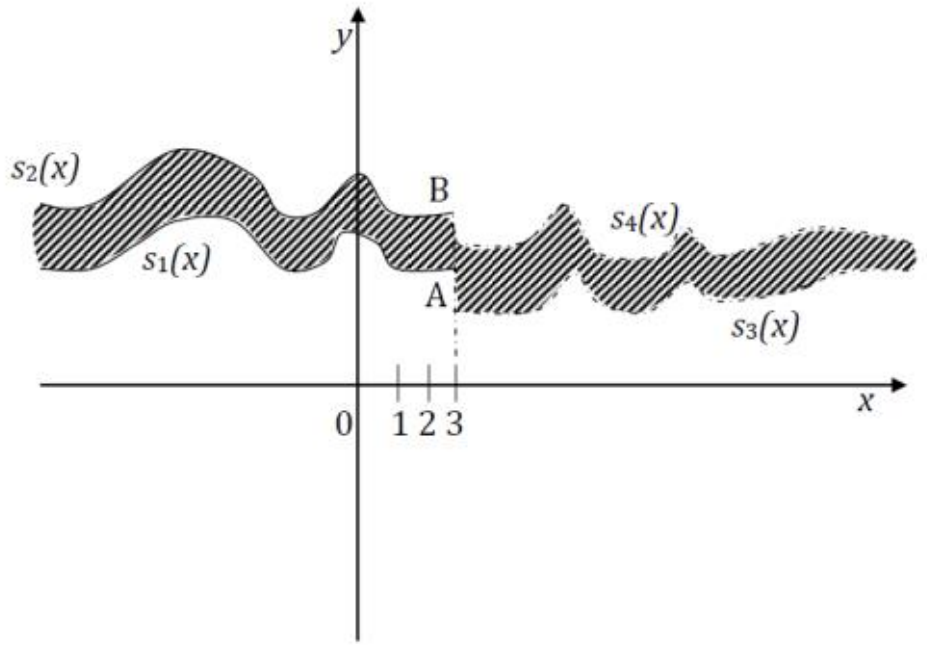

*Graph 7.*

For example, $s(3) = [s_1(3), s_2(3)]$, which is the vertical closed segment of line [AB].

In all above examples the indeterminacy occured into the values of function. But it is also possible to have indeterminacy into the argument of the function, or into both (the argument of the function, and the values of the function) as below.

8. Indeterminacy into the argument of the function:

$$r: \{1, 2, 3, 4\} \to \{5, 6, 7\} \quad (36)$$

$r(1) = 5, r(2) = 6,$

$r(3 \text{ or } 4) = 7$ {i. e. we do not know if $r(3) = 7$ or $r(4) = 4$}.

Another such example:

$$t: \{1, 2, 3, 4\} \to \{5, 6\} \quad (37)$$

$t(1) = 5,$ but $t(2 \text{ or } 3 \text{ or } 4) = 6.$

9. Indeterminacy into both:

$$\sqcup : \{1, 2, 3, 4\} \to \{5, 6, 7\} \quad (38)$$

$\sqcup\ (1 \text{ or } 2) = 5 \text{ or } 6 \text{ or } 7,$

which means that either U(1) = 5, or U(1) = 6, or U(1) = 7, or U(2) = 5, or U(2) = 6, or U(2) = 7;

$\sqcup$ (2 or 3 or 4) = 6 or 7.

Another example:

$v_1: \mathbb{R} \rightarrow \mathcal{P}(\mathbb{R}), v_1(x \text{ or } 2x) = 5x.$ (39)

Yet, this last neutrosophic function with indeterminacy into argument can be transformed, because $v_1(2x) = 5x$ is equivalent to $v_1(x) = 2.5x$ , into a neutrosophic function with indeterminacy into the values of the function only:

$v_2(x) = 2.5x \text{ or } 5x.$

Nor these last neutrosophic functions are relationships in a classical way.

# II.7. Discrete and Non-Discrete Indeterminacy

From another view point, there is a *discrete indeterminacy*, as for examples:

$f(2 \text{ or } 3) = 4$,
or $f(2) = 5 \text{ or } 6$,
or $f(2 \text{ or } 3) = 5 \text{ or } 6$;
and *non-discrete indeterminacy*, as for examples:

$f(7x \text{ or } 8x) = 63$,
or $f(x) = 10x^3$
or $20\sin(x)$,
or $f(x^2 \text{ or } 8x) = 16e^x$ and $\ln x$.

Depending on each type of indeterminacy we need to determine a specific *neutrosophic technic* in order to overcome that indeterminacy.

# II.8. Neutrosophic Vector-Valued Functions of Many Variables

We have given neutrosophic examples of real-valued functions of a real variable. But similar neutrosophic vector-valued functions of many variables there exist in any scientific space:

$$f: A_1 \times A_2 \times \dots \times A_n \to B_1 \times B_2 \times \dots B_m$$

$$f(x_1, x_2, \dots, x_n) = \begin{pmatrix} f_1(x_1, x_2, \dots, x_n), \\ f_2(x_1, x_2, \dots, x_n), \dots, \\ f_m(x_1, x_2, \dots, x_n) \end{pmatrix}. \qquad (40)$$

Sure $A_1, A_2, \dots, A_n$ and $B_1, B_2, \dots, B_n$ may be scientific spaces of any types.

Such neutrosophic vector-valued functions of many variables may have indeterminacy into their argument, into their values, or into both. And the indeterminacy can be discrete or non-discrete.

## II.9. Neutrosophic Implicit Functions

Similarly to the classical explicit and implicit function, there exist: *Neutrosophic Explicit Functions*, for example:

$$f(x) = x^2 \text{ or } x^2 + 1, \tag{41}$$

and *Neutrosophic Implicit Functions*, for example:

$$\{(x, y) \in \mathbb{R}^2 | e^x + e^y = 0 \text{ or } e^x + e^y = -1\}. \tag{42}$$

# II.10. Composition of Neutrosophic Functions

*Composition of Neutrosophic Functions* is an extension of classical composition of functions, but where the indeterminacy propagates.

For example:

$$f(x) = [\ln(x), \ln(3x)], \text{ for } x > 0, \tag{43}$$

and $$g(x) = \begin{cases} \frac{1}{x-5}, \text{if } x \neq 5; \\ 7 \text{ or } 9, \text{if } x = 5; \end{cases} \tag{44}$$

are both neutrosophic functions.

What is $(f \circ g)(5) =$?

$$(f \circ g)(5) = f\big(g(5)\big) = f(7 \text{ or } 9) =$$
$$[\ln 7, \ln 21] \text{ or } [\ln 9, \ln 27]. \tag{45}$$

Therefore, the discrete indeterminacy "7 or 9" together with the non-discrete (continous) indeterminacy " $[\ln(x), \ln(3x)]$ " have propagated into a double non-discrete (continuous) indeterminacy " $[\ln 7, \ln 21]$ or $[\ln 9, \ln 27]$ ".

But what is $(g \circ f)(5) =$?

$$(g \circ f)(5) = g\big(f(5)\big) = g([\ln 5, \ln 15]) =$$
$$\left[\frac{1}{\ln(15)-5}, \frac{1}{\ln(5)-5}\right] \approx [-0.43631, -0.29494]. \tag{46}$$

What is in general $(f \circ g)(x) =$?

$$(f \circ g)(x) = f\big(g(x)\big) = \begin{cases} f\left(\frac{1}{x-5}\right), \text{for } x \neq 5; \\ f(7 \text{ or } 9), \text{for } x = 5; \end{cases}$$
$$= \begin{cases} \left[\ln\left(\frac{1}{x-5}\right), \ln\left(\frac{3}{x-5}\right)\right], \text{for } x > 5; \\ \big[[\ln 7, \ln 21] \text{ or } [\ln 9, \ln 27]\big], \text{for } x = 5. \end{cases} \tag{47}$$

Since the domain of $f(\cdot)$ is $(0, \infty)$, one has $\frac{1}{x-5} > 0$, i.e. $x > 5$ for the first piecewise of $f \circ g$.

As we said before, a neutrosophic function $y = f(x)$ may have indeterminacy into its domain, or into its range, or into its relation between $x$ and $y$ (or into any two or three of them together).

## II.11. Inverse Neutrosophic Function

The *inverse of a neutrosophic function* is also a neutrosophic function, since the indeterminacy of the original neutrosophic function is transmitted to its inverse.

Example.

$$f(x) = \begin{cases} 2x + 1 \text{ or } 6x, \text{for } x \neq 0; \\ [1, 3], \text{for } x = 0; \end{cases} \quad (48)$$

or

$0 \neq x \longrightarrow 2x+1$ or $6x$;

$0 \longrightarrow [1, 3]$.

Let's find the inverse of the neutrosophic function $f(x)$.

$$y = 2x + 1 \text{ or } 6x, \text{for } x \neq 0. \quad (49)$$

Therefore $y = 2x + 1$ or $y = 6x$, for $x \neq 0$.

Interchange the variables: $x = 2y + 1$ or $x = 6y$, for $y \neq 0$.

Thus $x = 2y + 1$, whence $y = \frac{x-1}{2} \neq 0$, therefore $x \neq 1$, respectively: $x = 6y$, whence $y = \frac{x}{6} \neq 0$, therefore $x \neq 0$.

Hence, the inverse of the neutrosophic function $f(x)$ is:

$$f^{-1}(x) = \begin{cases} \frac{x-1}{2} \text{ or } \frac{x}{6}, & \text{for } x \neq 0 \text{ and } x \neq 1; \\ 0, & \text{for } x = [1, 3]. \end{cases} \quad (50)$$

Again, the inverse of a neutrosophic function:

$f = \mathbb{R} \to \mathbb{R}^2$

$f(x) = [2x + 1, 6x]$, for $x \in \mathbb{R}$,

or $x \to [2x + 1, 6x]$.

Simply, the inverse is:

$f^{-1}: \mathbb{R}^2 \to \mathbb{R}$

$f^{-1}([2x+1, 6x]) = x$, for all $x \in \mathbb{R}$,
or $[2x+1, 6x] \to x$. (51)

The inverse of the neutrosophic exponential function

$k(x) = 2^x$ or $x+1$

is $k^{-1}(x) = \log_2(x)$ or $\log_2(x+1)$. (52)

Similarly, the inverse of the neutrosophic logarithmic function

$h(x) = \log_{(0.09,\ 0.11)} x$

is $h^{-1}(x) = (0.09, 0.11)^x$. (53)

A classical function is invertible if and only if it is one-to-one (verifies the *Horizontal Line Test*).

Let's consider the classical function:

$f: \{1, 2, 3\} \to \{4, 5\}$ (54)

$f(1) = 4, f(2) = 5, f(3) = 5.$

This function is not one-to-one since it fails the horizontal line test at $y = 5$, since $f(2) = f(3)$. Therefore, this function is not classically invertible.

However, neutrosophically we can consider the neutrosophic inverse function

$f^{-1}(4) = 1, f^{-1}(5) = \{2, 3\}$,

$f^{-1}: \{4, 5\} \to \mathcal{P}(\{1, 2, 3\})$. (55)

For the graph of a neutrosophic inverse function $f^{-1}(x)$ we only need to reflect with respect to the symmetry axis $y = x$ the graph of the neutrosophic function $f(x)$.

The indeterminacy of a neutrosophic function is transmitted to its neutrosophic inverse function.

## Proposition

Any neutrosophic function is invertible.

*Proof.* If $f(x)$ fails the horizontal line test $f: A \to B$, *at* $y = b$, from the domain of definition of the neutrosophic function, we define the neutrosophic inverse function

$$f^{-1}(b) = \{a \in A, f(a) = b\}, f^{-1}: B \to A. \qquad (56)$$

Let $f: A \to B$ be a neutrosophic function. If the neutrosophic graph of $f$ contains the neutrosophic point $(C, D)$, where $C \subseteq A$ and $D \subseteq B$, then the graph of the neutrosophic inverse function $f^{-1}$ contains the neutrosophic point $(D, C)$.

A *neutrosophic point* is a generalization of the classical point $(c, d)$, where $c \in A$ an $d \in B$, whose dimension is zero. A neutrosophic point is in general a thick point, which may have the dimension 0, 1, 2 or more (depending on the space we work in).

As examples, $\alpha([1,2],[4,6])$ has dimension *2*:

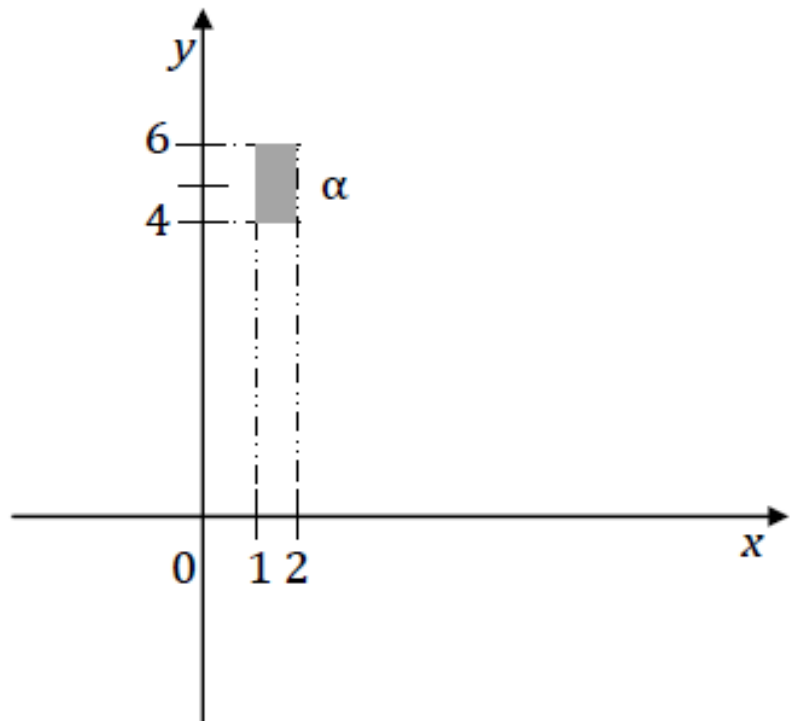

*Graph 8.*

or $\beta\big(3,(-1,1)\big)$ has the dimension 1:

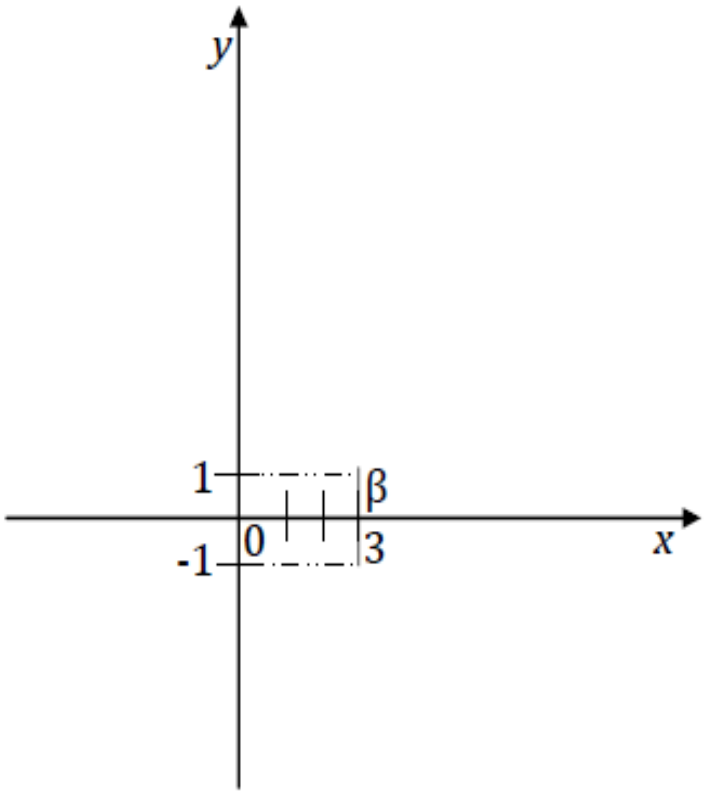

*Graph 9.*

or $\gamma(-2, \{-4, -3, -2\})$ has the dimension zero:

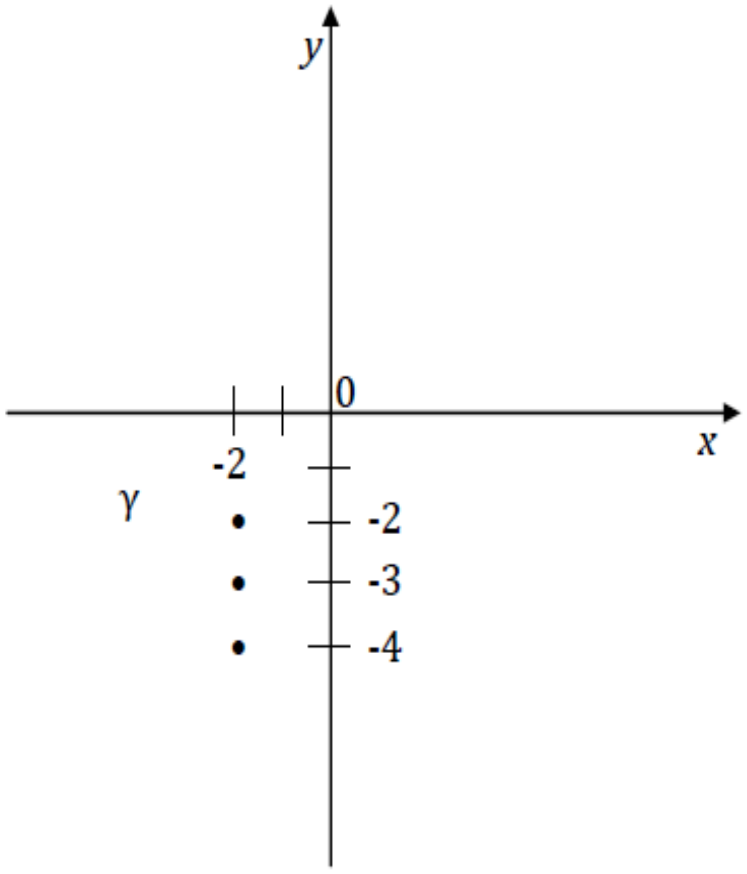

*Graph 10.*

while $\delta([2,3],[4,5],[0,4])$ has the dimension 3:

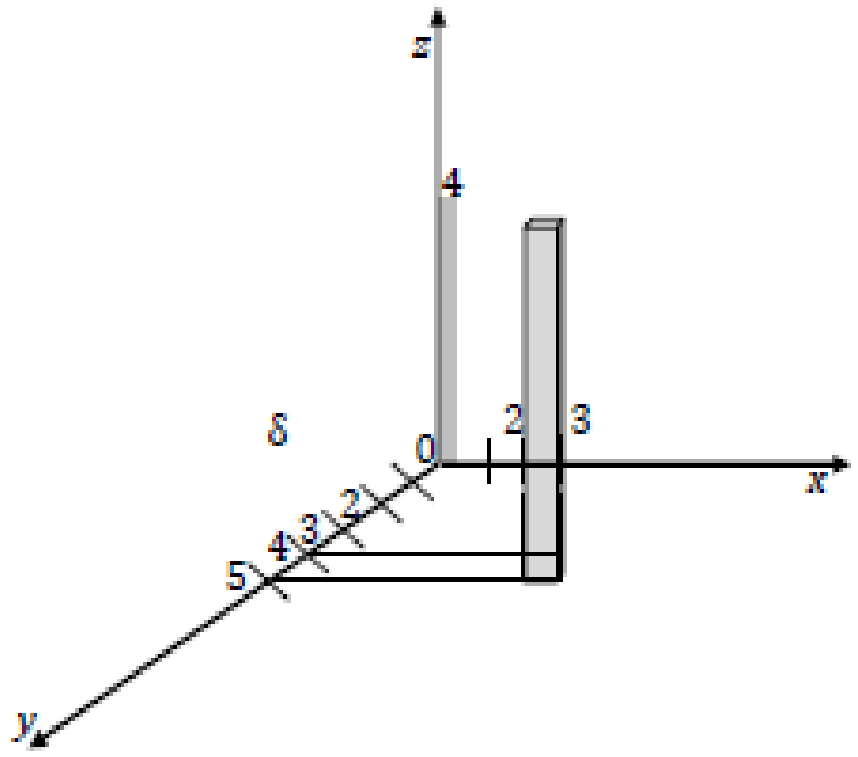

*Graph 11.*

# II.12. Zero of a Neutrosophic Function

Let $f: A \to B$. The zero of a neutrosophic function $f$ may be in general a set $S \subseteq A$ such $f(S) = 0$.

For example:

$f: \mathbb{R} \to \mathbb{R}$

$$f(x) = \begin{cases} x - 4, & x \notin [1,3] \\ 0, & x = [1,3] \end{cases}. \quad (57)$$

This function has a crisp zero, $x = 4$, since $f(4) = 4 - 4 = 0$, and an interval-zero $x = [1,3]$ since $f([1,3]) = 0$.

## II.13. Indeterminacies of a Function

By language abuse, one can say that any classical function is a neutrosophic function, if one considers that the classical function has a null indeterminacy.

# II.14. Neutrosophic Even Function

A *Neutrosophic Even Function:*

$$f: A \to B$$

has a similar definition to the classical even function:

$f(-x) = f(x)$, for all $x$ in $A$, (58)

with the extension that $f(-I) = f(I)$, where $I$ = indeterminacy.

For example:

$$f(x) = \begin{cases} x^2, & \text{for } x \notin \{-1, 1\}; \\ [0, 2], & \text{for } x = -1 \text{ or } 1. \end{cases} \quad (59)$$

Of course, for determinate

$x \in \mathbb{R} \setminus \{-1, 1\}, f(-x) = (-x)^2 = x^2 = f(x)$. (60)

While for the indeterminate $I = -1$ or $1$ one has

$-I = -(-1 \text{ or } 1) = 1 \text{ or } -1 = -1 \text{ or } 1$

whence $f(-I) = f(-1 \text{ or } 1) = [0, 2]$

and $f(I) = f(-1 \text{ or } 1) = [0, 2]$,

hence $f$ is a neutrosophic even function.

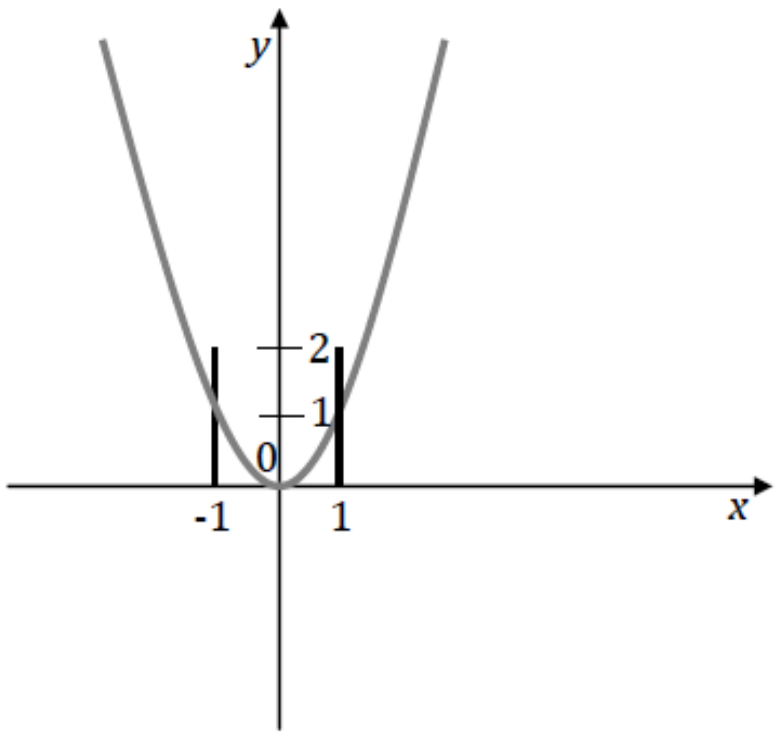

*Graph 12.*

As for classical even functions, the graph of a neutrosophic even function is symmetric, in a neutrosophic way, with respect to the y-axis, i.e. for a neutrosophic point P situated in the right side of the y-axis there exists a neutrosophic point P’ situated in the left side of the y-axis which is symmetric with P, and reciprocally.

We recall that the graph a neutrosophic function is formed by neutrosophic points, and a neutrosophic point may have not only the dimension 0 (zero), but also dimension 1, 2 and so on depending on the spaces the neutrosophic function is defined on and takes values in, and depending on the neutrosophic function itself.

# II.15. Neutrosophic Odd Function

Similarly, a *Neutrosophic Odd Function* $f: A \to B$ has a similar definition to the classical odd function:

$f(-x) = -f(x)$, for all $x$ in $A$, with the extension that $f(-I) = -f(I)$, where $I$ = indeterminacy.

For example:

$f: \mathbb{R} \to \mathbb{R}$

$$f(x) = \begin{cases} x \text{ and } x^3, & \text{for } x \neq 0; \\ -5 \text{ or } 5, & \text{for } x = 0. \end{cases} \quad (61)$$

The first piece of the function is actually formed by putting together two distinct functions.

Of course, for $x \neq 0$, $f(-x) = -x$, and $(-x)^3 = -x$, and $-x^3 = -(x \text{ and } x^3) = -f(x)$.

While for $x = 0$, one has:

$f(-0) = f(0) = -5 \text{ or } 5$;

$-f(0) = -(-5 \text{ or } 5) = 5 \text{ or } -5 = -5 \text{ or } 5$.

So, $f(-0) = -f(0)$, hence $f$ is a neutrosophic odd function.

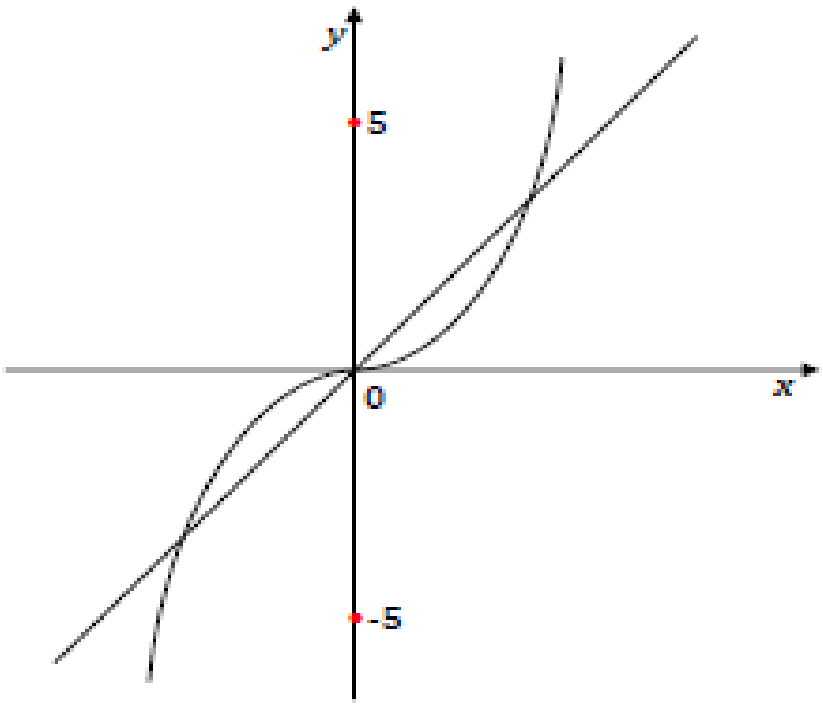

*Graph 13.*

Same thing: a neutrosophic odd function is neutrosophically symmetric with respect to the origin of the Cartesian system of coordinates.

# II.16. Neutrosophic Model

A model which has some indeterminacy is a neutrosophic model. When gathered data that describe the physical world is incomplete, ambiguous, contradictory, unclear, we are not able to construct an accurate classical model. We need to build an approximate (thick) model.

Using neutrosophic statistics, we plot the data and then design a neutrosophic regression method. The most common used such methods are the neutrosophic linear regression and the neutrosophic least squares regression.

For two neutrosophic variables, $x$ and $y$, representing the plotted data, one designs the best-fitting neutrosophic curve of the regression method. Instead of crisp data, as in classical regression, for example:

$$(x, y)\begin{Bmatrix} (1,2),(3,5),(4,8), \\ (-2,-4),(0,0),(-5,-11),\dots \end{Bmatrix}, \tag{62}$$

one works with set (approximate) data in neutrosophic regression:

$$(x, y) \in \begin{Bmatrix} (1,[2,2.2]),([2.5,3],5),\big([3.9,4),(8,8.1)\big), \\ (-2,-4),\big((0.0,0.1],(-0.1,0.0)\big), \\ (-5,(-10,-11)),\dots \end{Bmatrix} \tag{63}$$

and instead of obtaining, for example, a crisp linear regression as in classical statistics:

$$y = 2x - 1, \tag{64}$$

one gets a set-linear regression, for example:

$$y = [1.9,2]x - [0.9,1.1] \tag{65}$$

as in neutrosophic statistics.

# II.17. Neutrosophic Correlation Coefficient

The classical correlation coefficient $r$ is a crisp number between [-1, 1]. The neutrosophic correlation coefficient is a subset of the interval [-1, 1].

Similarly, if the subset of the neutrosophic correlation coefficient is more in the positive side of the interval [-1, 1], the neutrosophic variables $x$ and $y$ have a neutrosophic positive correlation, otherwise they have a neutrosphic negative correlation.

Of course, there is not a unique neutrosophic model to a real world problem. And thus, there are no exact neutrosophic rules to be employed in neutrosophic modelling. Each neutrosophic model is an approximation, and the approximations may be done from different points of view. A model might be considered better than others if it predicts better than others. But in most situations, a model could be better from a standpoint, and worse from another standpoint – since a real world problem normally depends on many (known and unknown) parameters.

Yet, a neutrosophic modelling of reality is needed in order to fastly analyse the alternatives and to find approximate optimal solutions.

# II.18. Neutrosophic Exponential Function

A *Neutrosophic Exponential Function* is an exponential function which has some indeterminacy [with respect to one or more of: its formula (base or exponent), or domain, or range].

If one has a classical exponential function

$g(x) = a^x$, with $a > 0$ and $a \neq 1$, (66)

then an indeterminacy with respect to the base can be, for example:

$f(x) = [0.9, 1.1]^x$, (67)

where "a" is an interval which even includes 1, and we get a thick function:

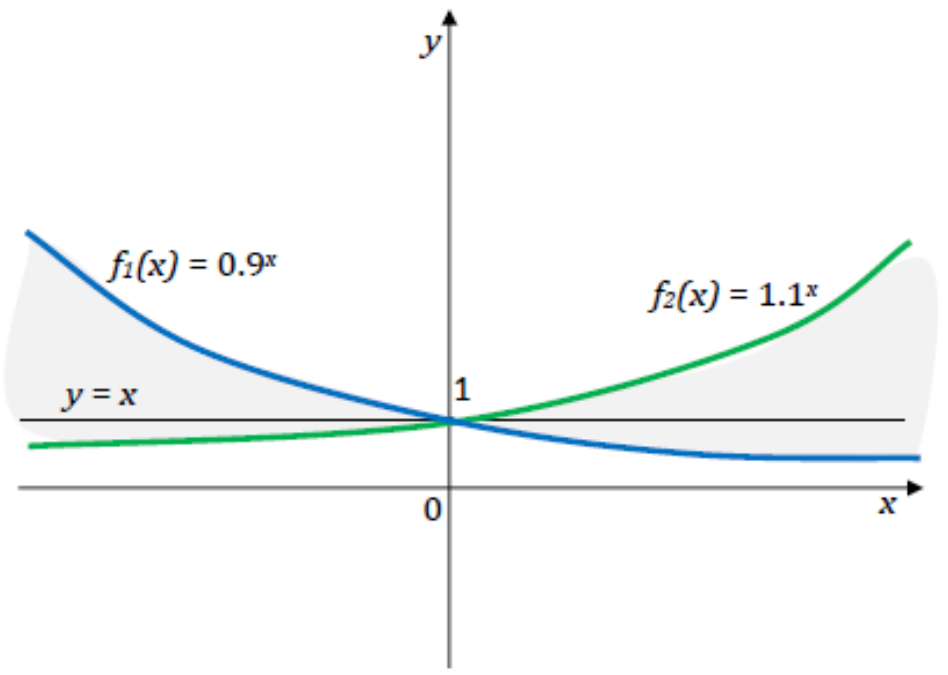

*Graph 14.*

or one may have indeterminacy with respect to the exponent:

$k(x) = 2^{x \text{ or } x+1}$. (68)

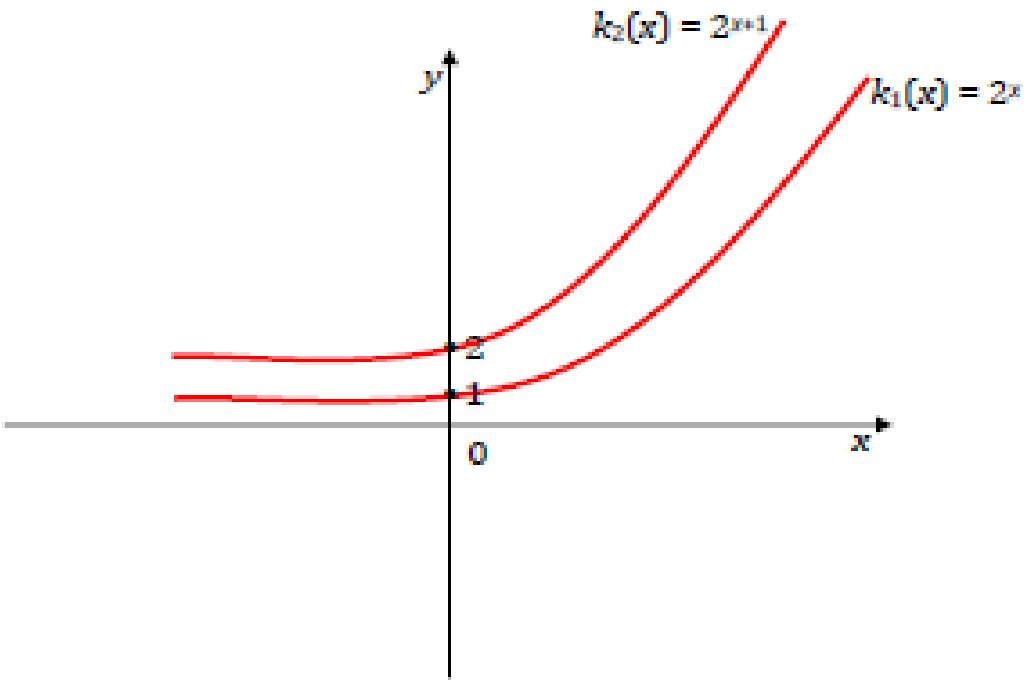

*Graph 15.*

For example: $k(1) = 2^{1 \text{ or } 1+1} = 2^1 \text{ or } 2^2 = 2 \text{ or } 4$ (we are not sure if it's 2 or 4). (69)

A third neutrosophic exponential function:

$l(x) = 2^{(x,\ x+1)}$ (70)

is different from $k(x)$ and has the graph:

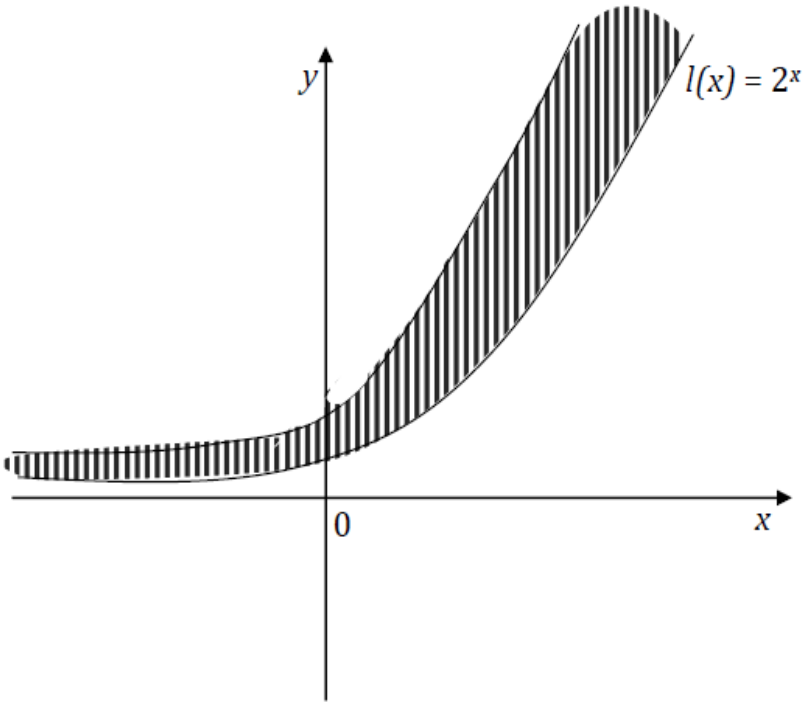

*Graph 16.*

which is a thick function. For example: $l(1) = 2^{(1,\ 1+1)} = 2^{(1,\ 2)} = (2^1, 2^2) = (2, 4)$, an open interval. (71)

# II.19. Neutrosophic Logarithmic Function

Similarly, a *Neutrosophic Logarithmic Function* is a logarithmic function that has some indeterminacy (with respect to one or more of: its formula, or domain, or range).

For examples:

$$f(x) = \log_{[2,3]} x = [\log_3 x, \log_3 x]. \tag{72}$$

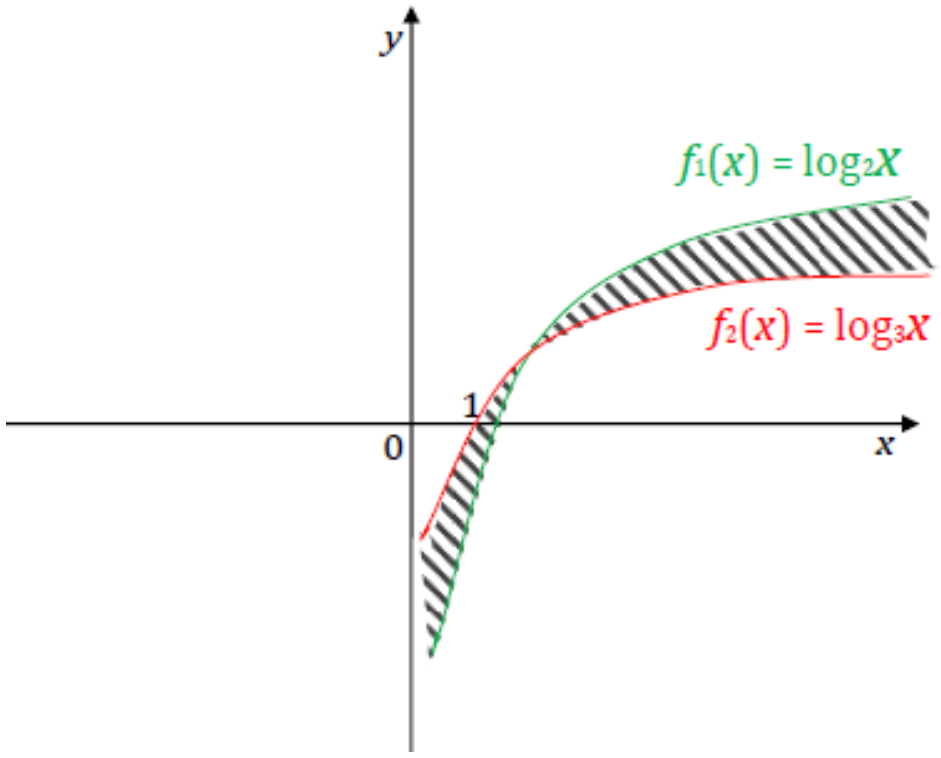

*Graph 17.*

or $g(x) = \ln(x, 2x) = \big(\ln(x), \ln(2x)\big)$ (73)

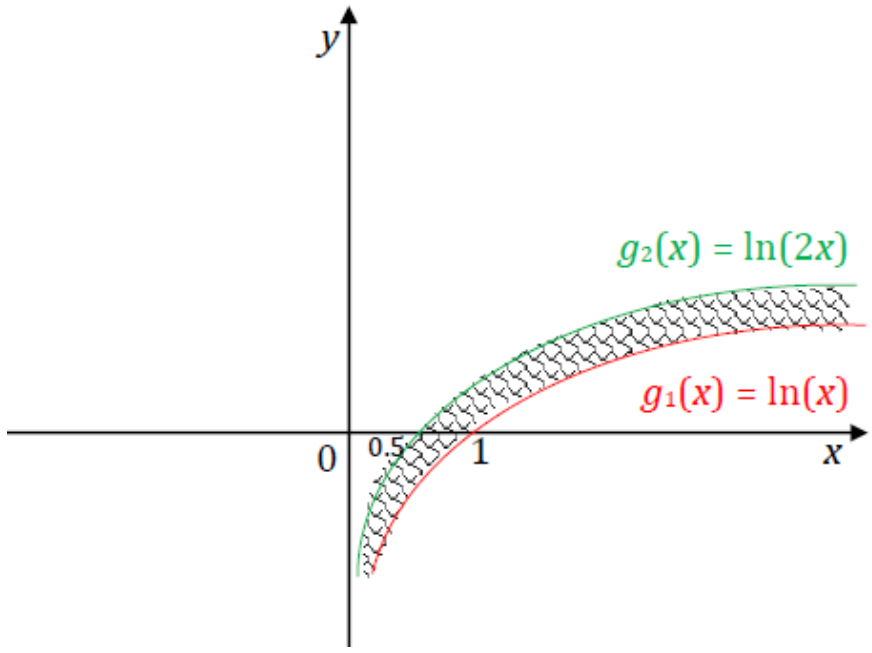

*Graph 18.*

or $h(x) = \log_{(0.09,11)} x$ (74)

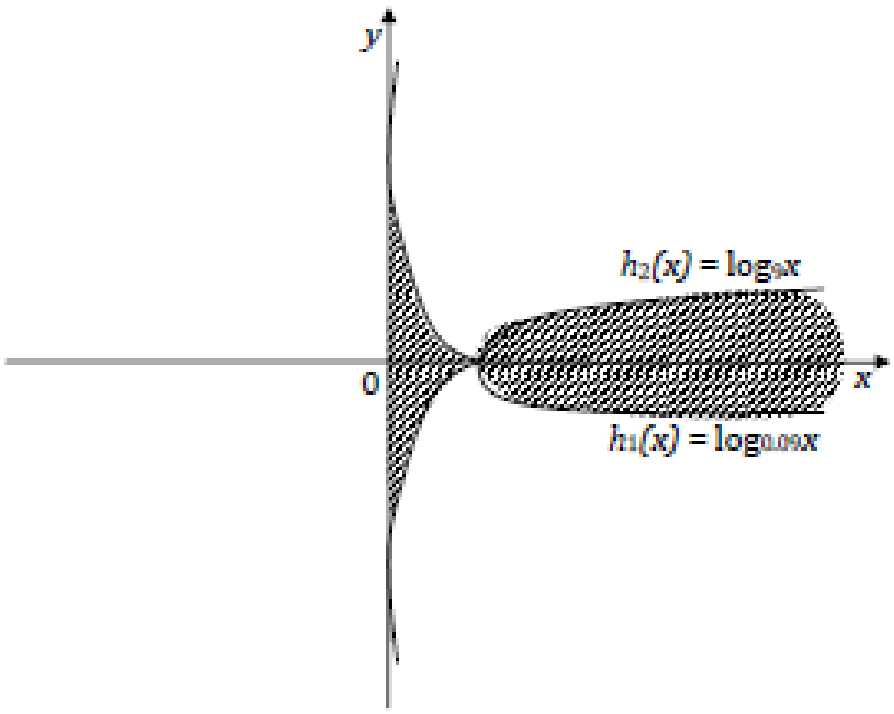

*Graph 19.*

## II.20. Composition of Neutrosophic Functions

In general, by composing two neutrosophic functions, the indeterminacy increases.

Example:

$$f_1(x) = x^3 \ or \ x^4$$
$$f_2(x) = [2.1, 2.5]^x$$

then

$$(f_1 \circ f_2)(x) = f_1\big(f_2(x)\big) = [2.1, 2.5]^{3x} \text{ or } [2.1, 2.5]^{4x}. \quad (75)$$

# III. Neutrosophic Calculus

# III.1. Neutrosophic Limit

*Neutrosophic Limit* means the limit of a neutrosophic function.

We extend the classical limit.

Let consider a neutrosophic function $f: \mathbb{R} \to \mathcal{P}(\mathbb{R})$ whose neutrosophic graph is below:

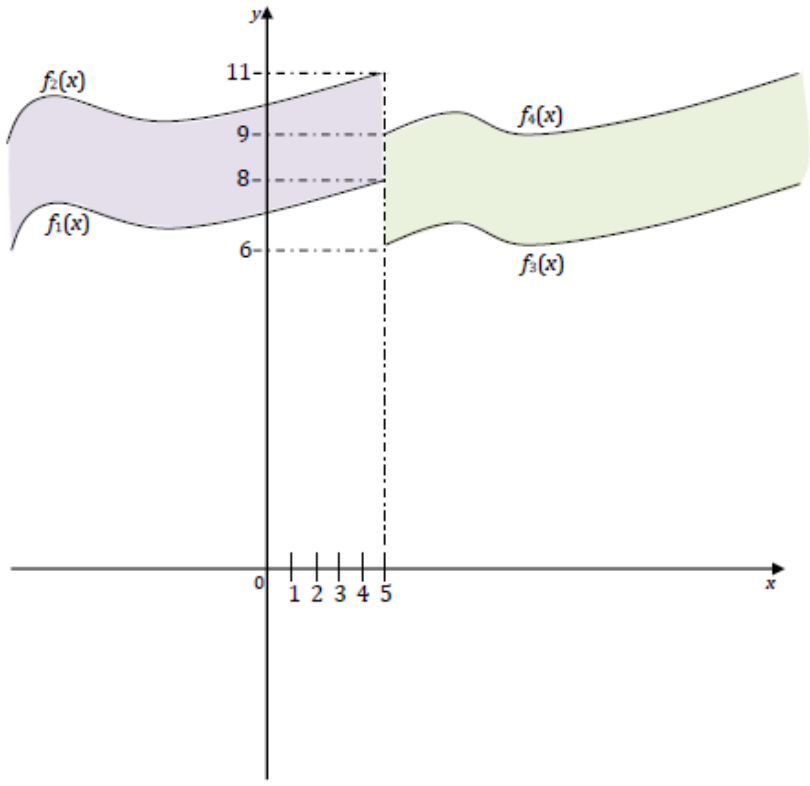

*Graph 20.*

$$f(x) = \begin{cases} [f_1(x), f_2(x)], \text{for } x \leq 5; \\ [f_3(x), f_4(x)], \text{for } x > 5, \end{cases} \quad (76)$$

is a neutrosophic piecewise-function.

Using the Neutrosophic Graphic Method, we get:

- The Neutrosophic Left Limit is

$$\lim_{\substack{x \to 5 \\ x < 5}} f(x) = [8, 11]; \quad (77)$$

- The Neutrosophic Right Limit is

$$\lim_{\substack{x \to 5 \\ x > 5}} f(x) = [6, 9]. \quad (78)$$

We introduce for the first time the notion of *neutrosophic mereo-limit*. Because the neutrosophic mereo-limit is the intersection of the neutrosophic left limit and the neutrosophic right limit [similarly as in the classical limit, where the left limit has to be equal to the right limit – which is equivalent to the fact that the intersection between the left limit (i.e. the set formed by a single finite number, or by $+\infty$, or by $-\infty$) and the right limit (i.e. also the set formed by a single finite number, or by $+\infty$, or by $-\infty$) is not empty], one has:

$$\lim_{x\to 5} f(x) = [8,11] \cap [6,9] = ]8,9]. \tag{79}$$

If the intersection between the neutrosophic left limit and the neutrosophic right limit is empty, then the neutrosophic mereo-limit does not exist.

Neutrosophic Limit of a function $f(x)$ does exist if the neutrosophic left limit coincides with the neutrosophic right limit. (We recall that in general the neutrosophic left and right limits are set, rather than numbers.) For example, the previous function does not have a neutrosophic limit since $[8,11] \not\equiv [6,9]$.

## Norm

We define a *norm*.

Let $\mu: \mathcal{P}(\mathbb{R}) \to \mathbb{R}^+$, where $\mathcal{P}(\mathbb{R})$ is the power set of $\mathbb{R}$, while $\mathbb{R}$ is the set of real numbers. (80)

For any set $\mathcal{S} \in \mathcal{P}(\mathbb{R})$,

$$\mu(\mathcal{S}) = \max\{|x|\}, x \in \mathcal{S} \cup Fr(\mathcal{S})\}, \tag{81}$$

where $|x|$ is the absolute value of $x$, and $Fr(\mathcal{S})$ is the frontier of $\mathcal{S}$,

or:

$$\mu(\mathcal{S}) = \max\{|inf\mathcal{S}|, |sup\mathcal{S}|\} \tag{82}$$

where $inf\mathcal{S}$ means the infimum of $\mathcal{S}$, and $sup\mathcal{S}$ means the supremum of $\mathcal{S}$.

Then:

$\mu(\mathcal{S}_1 + \mathcal{S}_2) = \max\{|inf\mathcal{S}_1 + inf\mathcal{S}_2|, |sup\mathcal{S}_1 + sup\mathcal{S}_2|\}$,

$$\mu(\alpha \cdot \mathcal{S}) = \max\{|\alpha| \cdot |inf\mathcal{S}|, |\alpha| \cdot |sup\mathcal{S}|\}, \quad (83)$$

where $\alpha \in \mathbb{R}$ is a scalar.

If the cardinality of the set $\mathcal{S}$ is 1, i.e. $\mathcal{S} = \{a\}, a \in \mathbb{R}$, then $\mu(\mathcal{S}) = \mu(a) = |a|$. (84)

We prove that $\mu(\cdot)$ is a *norm*.

$\mu: \mathcal{P}(\mathbb{R}) \to \mathbb{R}^+$,

$\forall \mathcal{S} \in \mathcal{P}(\mathbb{R}), \mu(\mathcal{S}) = \max\{|x|, x \in \mathcal{S} \cup Fr(\mathcal{S})\} = \max\{|inf\mathcal{S}|, |sup\mathcal{S}|\}$. (85)

$\mu(-\mathcal{S}) = \mu(-1 \cdot \mathcal{S}) = \max\{|-1| \cdot |inf\mathcal{S}|, |-1| \cdot |sup\mathcal{S}|\} = \max\{|inf\mathcal{S}|, |sup\mathcal{S}|\} = \mu(\mathcal{S})$. (86)

For a scalar $t$,

$\mu(t \cdot \mathcal{S}) = \max\{|t| \cdot |inf\mathcal{S}|, |t| \cdot |sup\mathcal{S}|\} = |t| \cdot \max\{|inf\mathcal{S}|, |sup\mathcal{S}|\} = |t| \cdot \mu(\mathcal{S})$. (87)

$\mu(S_1 + S_2) = max\{|infS_1 + infS_2|, |supS_1 + supS_2|\} \leq max\{|infS_1| + |infS_2|, |supS_1| + |supS_2|\} \leq max\{|infS_1|, |supS_1|\} + max\{|infS_2|, |supS_2|\} = \mu(S_1) + \mu(S_2)$. (88)

$\mu(S_1 - S_2) = \mu(S_1 + (-S_2)) \leq \mu(S_1) + \mu(-S_2) = \mu(S_1) + \mu(S_2)$. (89)

# III.2. Appropriateness Partial-Distance (Partial-Metric)

Let A and B be two sets included in $\mathbb{R}$, such that $inf A$, $supA$, $inf B$, and $supB$ are finite numbers.

Then the appropriate partial-distance (partial-metric) between A and B is defined as:

$\eta : \mathbb{R}^2 \longrightarrow \mathbb{R}^+$

$\eta(A, B) = max\{|infA\text{-}infB|, |supA\text{-}supB|\}.)$ (90)

In other words, the appropriateness partial-distance measures how close the *inf's* and *sup's* of two sets (i.e. the two sets corresponding extremities) are to each other.

# III.3. Properties of the Appropriateness Partial-Distance

For any A, B, C $\subset \mathbb{R}$, such that $inf A$, $sup A$, $inf B$, $sup B$, $inf C$, and $sup C$ are finite numbers, one has:

a) $\eta(A, B) \geq 0.$ (91)

b) $\eta(A, A) = 0.$ (92)

But if $\eta(A, B) = 0$ it does not result that A ≡ B, it results that $infA = infB$ and $supA = supB$.

For example, if $A = \{3, 4, 5, 7\}$ and $B = (3, 7]$, then $infA = infB = 3$ and $supA = supB = 7$, whence $\eta(A, B) = 0$, but $A \not\equiv B$. (93)

Therefore, this distance axiom is verified only partially by $\eta$.

c) $\eta(A, B) = \eta(B, A).$ (94)

d) $\eta(A, B) \leq \eta(B, C) + \eta(C, A).$ (95)

Proof of d):

$\eta(A, B) = max\{|infA-infB|, |supA-supB|\}$

$= max\{|infA-infC + infC -infB|, |supA-supC+supC-supB|\}.$ (96)

But $|infA-infC + infC -infB| \leq |infA-infC| + |infC -infB| = |infB-infC| + |infC -infA|$ (97)

and similarly

$|supA-supC+supC-supB| \leq |supA-supC|+|supC-supB| = |supB-supC|+|supC-supA|$ (98)

whence

$max\{|infA-infC + infC-infB|, |supA-supC+supC-supB|\} \leq max\{|infB-infC|, |supB-supC|\} + max\{|infC-infA|, |supC-supA|\} = \eta(B, C)+ \eta(C, A).$ (99)

*e)* If $A = \{a\}$ and $B = \{b\}$, with $a, b \in \mathbb{R}$, i.e. A and B contain only one element each, then:
$\eta(A, B) = |a-b|$. (100)

f) If $A$ and $B$ are real (open, closed, or semi-open/semi-closed) intervals, $A = [a_1, a_2]$ and $B = [b_1, b_2]$, with $a_1 < a_2$ and $b_1 < b_2$, then $\eta(A, B) = max\{|a_1-b_1|, |a_2-b_2|\}$. (101)

## III.4. Partial-Metric Space

Let's have in general:

$\eta: \mathcal{M} \to R^+$, where $\mathcal{M}$ is a non-empty set.

The function $\eta$ is a *partial-metric* (partial-distance) on $\mathcal{M}$,

$$\eta(A, B) = max\{|infA - infB|, |supA - supB|\} \quad (102)$$

and the space $\mathcal{M}$ endowed with $\eta$ is called a *partial-metric space*.

This partial-metric space $\eta$ is a generalization of the metric $d$, defined in interval analysis:

$d: S \to S$, where $S$ is any real set, and

$$d([a, b], [c, d]) = max\{|a - c|, |b - d|\}, \quad (103)$$

with $a \leq b$ and $c \leq d$, because $\eta$ deals with all kinds of sets, not only with intervals as in integer analysis.

Remarkably,

$$\eta(A, 0) = max\{|infA - 0|, sup|A - 0|\} =$$

$$max\{|infA|, |supA|\} = \mu(A), \quad (104)$$

which is the norm of $A$.

# III.5. ε − δ Definition of the Neutrosophic Left Limit

Let $f$ be a neutrosophic function, $f: P(\mathbb{R}) \longrightarrow P(\mathbb{R})$.

*The ε − δ definition of the Neutrosophic Left Limit* is an extension of the classical left limit definition, where the absolute value |·| is replace by $\eta(\cdot)$. Also, instead of working with scalars only, we work with sets (where a "set" is view as an approximation of a "scalar").

Therefore,

$$\lim_{\substack{x \to c \\ x < c}} f(x) = L \tag{105}$$

is equivalent to $\forall \varepsilon > 0$, $\exists \delta = \delta(\varepsilon) > 0$, such that if $\eta(x, c)_{x<c} < \delta$, then $\eta(f(x), L)_{x<c} < \varepsilon$.
(106)

*The ε − δ definition of the Neutrosophic Right Limit*.

$$\lim_{\substack{x \to c \\ x > c}} f(x) = L \tag{107}$$

is equivalent to $\forall \varepsilon > 0$, $\exists \delta = \delta(\varepsilon) > 0$, such that if $\eta(x, c)_{x>c} < \delta$, then $\eta(f(x), L)_{x>c} < \varepsilon$.
(108)

And, in general, the *ε − δ definition of the Neutrosophic Limit.*

$$\lim_{x \to c} f(x) = L$$

is equivalent to $\forall \varepsilon > 0$, $\exists \delta = \delta(\varepsilon) > 0$, such that if $\eta(x, c) < \delta$, then $\eta(f(x), L) < \varepsilon$.
(109)

# III.6. Example of Calculating the Neutrosophic Limit

In our *previous example*, with $c = 5$, let $\varepsilon > 0$, then

$$\eta([f_1(x), f_2(x)], [8, 11]) = \max_{\substack{\eta(x-5)<\delta \\ x<5}} \{|inf[f_1(x), f_2(x)] - inf[8, 11]|, |sup[f_1(x), f_2(x)] - sup[8, 11]|\} = \max_{\substack{\eta(x-5)<\delta \\ x<5}} \{|f_1(x) - 8|, |f_2(x) - 11|\} < \varepsilon. \quad (110)$$

$\eta(x, 5) < \delta$ means $|x - 5| < \delta$ as in classical calculus.

$$\max_{\substack{\eta(x-5)<\delta \\ x<5}} \{|f_1(x) - 8|, |f_2(x) - 11|\} < \varepsilon$$

means $|f_1(x) - 8| < \varepsilon$, and $|f_2(x) - 11| < \varepsilon$, when $|x - 5| < \delta$ and $x \leq 5$. (111)

# III.7. Particular Case of Calculating the Neutrosophic Limit

Suppose, as a *particular case* of the previous example, that $f_1(x), f_2(x), f_3(x), f_4(x)$ are piecewise functions, such that in a left or right neighborhood of $x = 5$ they are:

$f_1(x) = -x^2 + 6x + 3$, for $x \in [4, 5]$; (112)

$f_2(x) = x^3 - 114$, for $x \in [4, 5]$; (113)

$f_3(x) = x + 1$, for $x \in [5, 6]$; (114)

$f_4(x) = 3x - 6$, for $x \in [5, 6]$. (115)

Therefore,

$|f_1(x) - 8| = |-x^2 + 6x + 3 - 8| = |-(x - 5)(x - 1)| = |(x - 5)(x - 1)| < \frac{\varepsilon}{4}(4) = \varepsilon$; we take $\delta = \frac{\varepsilon}{4}$, because $x - 1 \leq 4$, since $x \in [4, 5]$. (116)

And $|f_2(x) - 11| = |x^3 - 114 - 11| = |(x - 5)(x^2 + 5x + 25)| < \frac{\varepsilon}{75}(75) = \varepsilon$; we take $\delta = \frac{\varepsilon}{75}$, because $x^2 + 5x + 25 \leq (5)^2 + 5(5) + 25 = 75$, since $x \in [4, 5]$. (117)

We got that for any $\varepsilon > 0$, there exists $\delta = min\left\{\frac{\varepsilon}{4}, \frac{\varepsilon}{75}\right\} = \frac{\varepsilon}{75}$. Whence it results the *neutrosophic left limit*.

Similarly for the *neutrosophic right limit* in this example.

Let $\varepsilon > 0$. Then

$$\eta([f_3(x), f_4(x)], [6, 9]]) = \max_{\substack{\eta(x-5)<\delta \\ x>5}} \{|inf[f_3(x), f_4(x)] - \inf[6, 9]|, |sup[f_3(x), f_4(x)] - \sup[6, 9]|\} = \max_{\substack{\eta(x-5)<\delta \\ x>5}} \{|f_3(x) - 6|, |f_4(x) - 9|\} < \varepsilon, \quad (118)$$

which means

$|f_3(x) - 6| < \varepsilon$, and $|f_4(x) - 9| < \varepsilon$,
when $|x - 5| < \delta$ and $x > 5$.

Therefore:

$$|f_3(x) - 6| = |x + 1 - 6| = |x - 5| < \frac{\varepsilon}{1}(1) = \varepsilon;$$

we take $\delta = \frac{\varepsilon}{1} = \varepsilon$. (119)

And:

$$|f_4(x) - 9| = |3x - 6 - 9| = |3(x - 5)| < \frac{\varepsilon}{3} \cdot (3) = \varepsilon;$$

we take $\delta = \frac{\varepsilon}{3}$. (120)

We got that for any $\varepsilon > 0$, there exists

$$\delta = min\left\{\varepsilon, \frac{\varepsilon}{3}\right\} = \frac{\varepsilon}{3}, \quad (121)$$

whence it results the neutrosophic right limit.

Then we intersect the neutrosophic left and right limits to get the neutrosophic mereo-limit. We observe that the neutrosophic limit does not exist of this function, since if we take $\varepsilon = 0.1 > 0$, there exist no $\delta = \delta(\varepsilon) > 0$ such that if $|x - 5| < \delta$ to get

$$\eta([f_1(x), f_2(x)], [8, 9]) < 0.1 \quad (122)$$

not even

$$\eta([f_3(x), f_4(x)], [8, 9]) < 0.1 \quad (123)$$

since in tiny neighborhood of 5 the absolute values of differences $|f_2(x) - 9|$ and $|f_3(x) - 8|$ are greater than 1.

# III.8. Computing a Neutrosophic Limit Analytically

Let's consider the below limit:

$$\lim_{x\to -3} \frac{x^2+3x-[1,2]x-[3,6]}{x+3} \tag{124}$$

We substitute $x$ for -3, and we get:

$$\lim_{x\to -3} \frac{(-3)^2+3\cdot(-3)-[1,2]\cdot(-3)-[3,6]}{-3+3} = \frac{9-9-[1\cdot(-3),2\cdot(-3)]-[3,6]}{0} = \frac{0-[-6,-3]-[3,6]}{0} = \frac{[3,6]-[3,6]}{0} = \frac{[3-6,6-3]}{0} = \frac{[-3,3]}{0}, \tag{125}$$

which has un undefined operation $\frac{0}{0}$, since $0 \in [-3,3]$.

Then we factor out the numerator, and simplify:

$$\lim_{x\to -3} \frac{x^2+3x-[1,2]x-[3,6]}{x+3} = \lim_{x\to -3} \frac{(x-[1,2])\cdot(x+3)}{(x+3)} = \lim_{x\to -3}(x-[1,2]) = -3-[1,2] = [-3,-3]-[1,2] = -([3,3]+[1,2]) = [-5,-4]. \tag{126}$$

We can check the result considering classical crisp coefficients instead of interval-valued coefficients.

For examples:

a) Taking the infimum of the intervals [1,2] and respectively [3,6], i.e. 1 and respectively 3, we have:

$\lim_{x\to -3}\frac{x^2+3x-1x-3}{x+3} =$

$\lim_{x\to -3}\frac{x^2+2x-3}{x+3} = \lim_{x\to -3}\frac{(x+3)(x-1)}{x+3} = \lim_{x\to -3}(x-1)$ =-3-1 = -4 $\in[-5,-4]$. (127)

b) Taking the supremum of the intervals [1,2] and respectively [3,6], i.e. 2 and respectively 6, we have:

$\lim_{x\to -3}\frac{x^2+3x-2x-6}{x+3} =$

$\lim_{x\to -3}\frac{x^2+x-6}{x+3} = \lim_{x\to -3}\frac{(x+3)(x-2)}{x+3} = \lim_{x\to -3}(x-2)$ =-3-2 = = -5 $\in[-5,-4]$. (128)

c) Taking the midpoints of the intervals [1,2] and respectively [3,6], i.e. 1.5 and respectively 4.5, we have:

$\lim_{x\to -3}\frac{x^2+3x-1.5x-4.5}{x+3} =$

$\lim_{x\to -3}\frac{x^2+1.5x-4.5}{x+3} = \lim_{x\to -3}\frac{(x+3)(x-1.5)}{x+3} = \lim_{x\to -3}(x-1.5) =$ -3-1.5 = -4.5 $\in[-5,-4]$. (129)

d) In general, taking α ∈ [1,2] and respectively 3α ∈ [3,6], one has:

$\lim_{x \to -3} \frac{x^2+3x-\alpha x-3\alpha}{x+3} =$
$\lim_{x \to -3} \frac{x^2+(3-\alpha)x-3\alpha}{x+3} = \lim_{x \to -3} \frac{(x+3)(x-\alpha)}{x+3} = \lim_{x \to -3} (x -$
$\alpha)$ = -3- α ∈ [-3,-3]-[1,2]  { since α ∈ [1,2] }
= [-3-2, -3-1] = [-5, -4]. (130)

So, we got the same result.

# III.9. Calculating a Neutrosophic Limit Using the Rationalizing Technique

$$\lim_{x\to 0}\frac{\sqrt{(4,5)\cdot x+1}-1}{x}=\frac{\sqrt{(4,5)\cdot 0+1}-1}{0}$$
$$=\frac{\sqrt{[4\cdot 0,5\cdot 0]+1}-1}{0}$$
$$=\frac{\sqrt{[0,0]+1}-1}{0}=\frac{\sqrt{0+1}-1}{0}=\frac{0}{0}$$
$$=\text{undefined}. \quad (131)$$

Multiply with the conjugate of the numerator:

$$\lim_{x\to 0}\frac{\sqrt{[4,5]x+1}-1}{x}\cdot\frac{\sqrt{[4,5]x+1}+1}{\sqrt{[4,5]x+1}+1}$$
$$=\lim_{x\to 0}\frac{\left(\sqrt{[4,5]x+1}\right)^2-(1)^2}{x\left(\sqrt{[4,5]x+1}+1\right)}$$
$$=\lim_{x\to 0}\frac{[4,5]\cdot x\cancel{+1-1}}{x\cdot\left(\sqrt{[4,5]x+1}+1\right)}$$
$$=\lim_{x\to 0}\frac{[4,5]\cancel{\cdot x}}{\cancel{x\cdot}\left(\sqrt{[4,5]x+1}+1\right)}$$
$$=\lim_{x\to 0}\frac{[4,5]}{\left(\sqrt{[4,5]x+1}+1\right)}$$
$$=\frac{[4,5]}{\left(\sqrt{[4,5]\cdot 0+1}+1\right)}=\frac{[4,5]}{\sqrt{1}+1}$$
$$=\frac{[4,5]}{2}=\left[\frac{4}{2},\frac{5}{2}\right]=[2,2.5]. \quad (132)$$

Similarly we can check this limit in a classical way considering a parameter $\alpha \in [4,5]$ and computing the limit by multiplying with the conjugate of the numerator:

$$\lim_{x \to 0} \frac{\sqrt{\alpha \cdot x + 1} - 1}{x} = \frac{\alpha}{2} \in [4,5]/2 = [2, 2.5]. \quad (133)$$

# III.10. Neutrosophic Mereo-Continuity

We now introduce for the first time the notion of *neutrosophic mereo-continuity. A neutrosophic function* $f(x)$ is *mereo-continuous at a given point* $x = c$, where

$f: A \to B$

if the intersection of the neutrosophic left limit, neutrosophic right limit, and $f(c)$ is nonempty:

$$\left\{\lim_{\substack{x\to c\\ x<c}} f(x)\right\} \cap \left\{\lim_{\substack{x\to c\\ x>c}} f(x)\right\} \cap \{f(c)\} \neq 0. \quad (134)$$

A *neutrosophic function* $f(x)$ is *mereo-continuous on a given interval* $[a, b]$, if there exist the classical points $A \in \{f(a)\}$ and $B \in \{f(b)\}$ that can be connected by a continuous classical curve which is inside of $f(x)$.

Also, the classical definition can be extended in the following way: A *neutrosophic function* $f(x)$ is *mereo-continuous on a given interval* $[a, b]$, if $f(x)$ is neutrosophically continuous at each point of $[a, b]$.

A *neutrosophic function* $f(x)$ is *continuous* at a given point $x = c$ if:

$$\lim_{\substack{x\to c\\ x>c}} f(x) \equiv \lim_{\substack{x\to c\\ x<c}} f(x) \equiv f(c). \quad (135)$$

We see that the previous neutrosophic function is mereo-continuous at $x = 5$ because:

$$\left\{\lim_{\substack{x\to 5\\ x<5}} f(x)\right\} \cap \left\{\lim_{\substack{x\to 5\\ x>5}} f(x)\right\} \cap \{f(5)\} = [8, 11] \cap$$

$$[6, 9] \cap [8, 11] = [8, 9] \neq \phi. \quad (136)$$

# III.11. Neutrosophic Continuous Function

A neutrosophic function $f: \mathcal{M}_1 \to \mathcal{M}_2$ is *continuous* at a *neutrosophic point* $x = c$ if:

$\forall \varepsilon > 0, \exists\, \delta = \delta(\varepsilon) > 0,$ (137)

such that for any $x \in \mathcal{M}_1$ such that $\eta(x, c) < \delta$ one has $\eta\big(f(x), f(c)\big) < \varepsilon.$ (138)

(We recall that a "neutrosophic point" $x = c$ is in general a set $c \in \mathcal{M}_1$, while $\mathcal{M}_1$ and $\mathcal{M}_2$ are sets of sets.)

# III.12. Neutrosophic Intermediate Value Theorem

Let $f: A \to P(A), f(x) = [a_x, b_x] \subseteq A$, where $[a_x, b_x]$ is an interval. (139)

Let

$inf\{f(a)\} = a_1;$

$sup\{f(a)\} = a_2;$

$inf\{f(b)\} = b_1;$

$sup\{f(b)\} = b_2.$

Suppose $min\{a_1, a_2, b_1, b_2\} = m$,

and $max\{a_1, a_2, b_1, b_2\} = M$.

If $f(x)$ is a neutrosophic mereo-continuous function on the closed interval $[a, b]$, and $k$ is a number between $m$ and $M$, with $m \neq M$, then there exists a number $c \in [a, b]$ such that: $\{f(c)\} \ni k$ (i.e. the set $\{f(c)\}$ contains $k$), or $k \in \{f(c)\}$.

An extended version of this theorem is the following:

If $f(x)$ is a neutrosophic mereo-continuous function of the closed interval [a, b], and $\langle k_1, k_2 \rangle$ is an interval included in the interval $[m, M]$, with $m \neq M$, then there exist $c_1, c_2, \dots, c_m$ in $[a, b]$, where $m \geq 1$, such that $\langle k_1, k_2 \rangle \subseteq f(c_1) \cup f(c_2) \cup \dots \cup f(c_m)$.

Where by $\langle \alpha, \beta \rangle$ we mean any kind of closed, open or half-closed and half-open intervals: $[\alpha, \beta]$, or $(\alpha, \beta)$, or $[\alpha, \beta)$, or $(\alpha, \beta]$.

# III.13. Example for the Neutrosophic Intermediate Value Theorem

Let $g(x) = [g_1(x), g_2(x)]$, where $g: \mathbb{R} \to \mathbb{R}^2$, and $g_1, g_2: \mathbb{R} \to \mathbb{R}$.

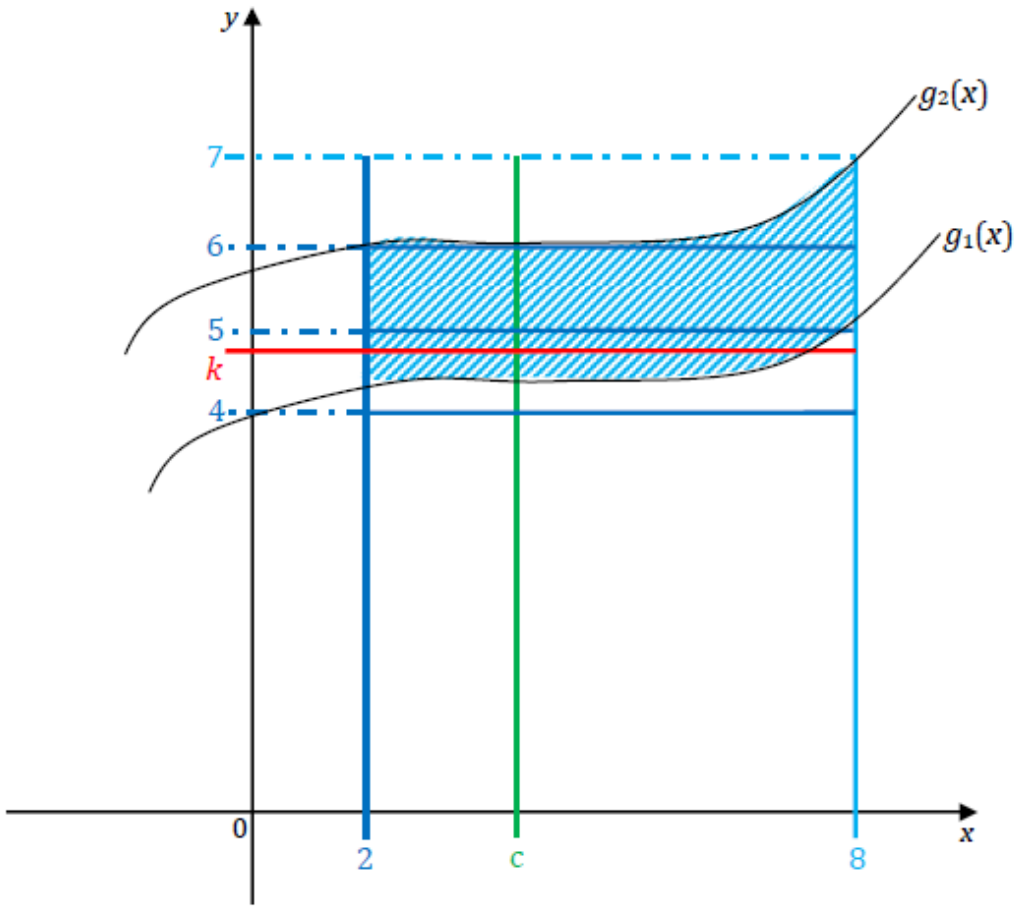

*Graph 21.*

$g$ is neutrosophically continuous on the interval $[2, 8]$.

Let $m = min\{4, 5, 6, 7\} = 4$,

and $M = max\{4, 5, 6, 7\} = 7$, and let $k \in [4, 7]$.

Then there exist many values of $c \in [2, 8]$ such that $\{g(c)\} \ni k$. See the green vertical line above, $x = c$. For example $c = 4 \in [2, 8]$. The idea is that if $k \in [4, 7]$ and we draw a horizontal red line $g = k$, this horizontal red line will intersect the shaded blue area which actually represents the neutrosophic graph of the function $g$ on the interval $[2, 8]$.

# III.14. Example for the Extended Intermediate Value Theorem

Let $h(x) = [h_1(x), h_2(x)]$, where $h: \mathbb{R} \to \mathbb{R}^2$, and $h_1, h_2: \mathbb{R} \to \mathbb{R}$. $h$ is neutrosophically continuous on the interval $[3, 12]$.

Let $m = min\{6, 8, 10, 12.5\} = 6$,
and $M = max\{6, 8, 10, 12.5\} = 12.5$,
and let $[k_1, k_2] \in [6.5, 12] \subset [6, 12.5]$.

Then there exist $c_1 = 8 \in [3, 12]$ and $c_2 = 10 \in [3, 12]$ such that

$$h(c_1) \cup h(c_2) = h(8) \cup h(10) = [6.5, 11] \cup [9.5, 12] = [6.5, 12] = [k_1, k_2]. \quad (140)$$

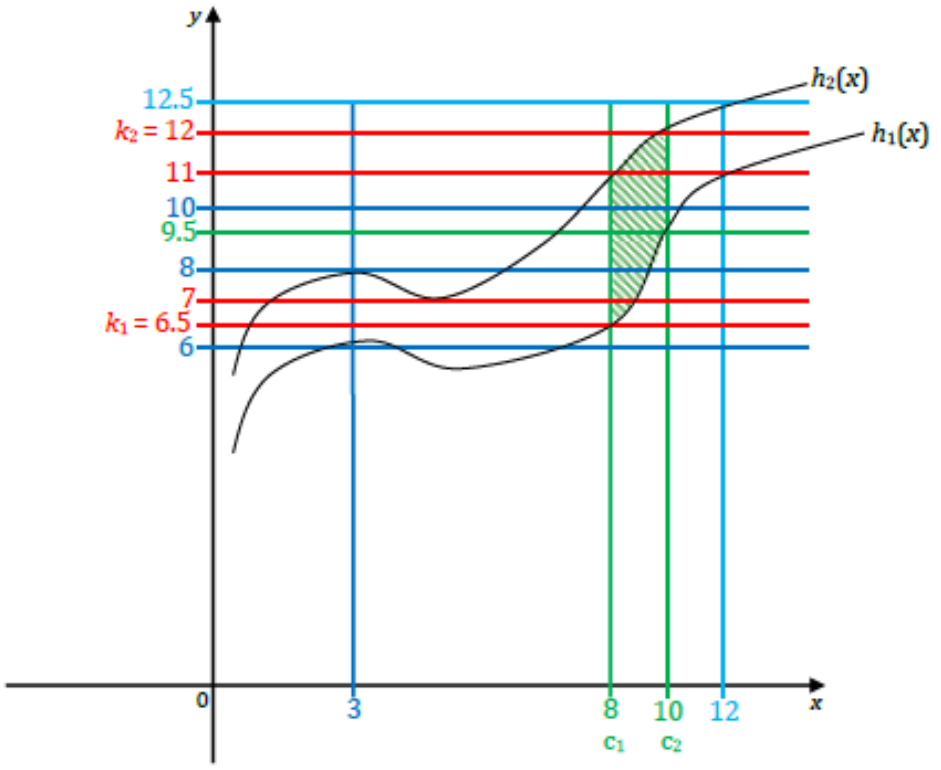

*Graph 22.*

## Remark

The more complicated (indeterminate) is a neutrosophic function, the more complex the neutrosophic intermediate value theorem becomes.

Actually, for each class of neutrosophic function, the neutrosophic intermediate value theorem has a special form.

As a General Remark, we have:

For each class of neutrosophic functions a theorem will have a special form.

## III.15. Properties of Neutrosophic Mereo-Continuity

1. A neutrosophic $f(x)$ is mereo-continuous on the interval $[a,b]$, if it's possible to connect a point of the set $\{f(a)\}$ with a point of the set $\{f(b)\}$ by a continuous classical curve $\mathbb{C}$ which is included in the (thick) neutrosophic function $f(x)$ on the interval $[a,b]$.

2. If $\alpha \neq 0$ is a real number, and $f$ is a neutrosophic mereo-continuous function at $x = c$, then $\alpha \cdot f$ is also a neutrosophic mereo-continuous function at $x = c$.

### Proof

$$\lim_{\substack{x\to c\\ x<c}}[a \cdot f(x)] \cap \lim_{\substack{x\to c\\ x>c}}[a \cdot f(x)] \cap \{\alpha \cdot f(c)\} =$$

$$\left\{\alpha \cdot \lim_{\substack{x\to c\\ x<c}}[f(x)]\right\} \cap \left\{\alpha \cdot \lim_{\substack{x\to c\\ x>c}}[f(x)]\right\} \cap \{\alpha \cdot f(c)\} = \alpha \cdot$$

$$\left(\left\{\lim_{x\to c}[f(x)]\right\} \cap \left\{\lim_{\substack{x\to c\\ x>c}}[f(x)]\right\} \cap \{f(c)\}\right) \neq \emptyset, \qquad (141)$$

because $\alpha \neq 0$, and $\left\{\lim_{\substack{x\to c\\ x<c}}[f(x)]\right\} \cap \lim_{\substack{x\to c\\ x>c}}[f(x)] \cap \{f(c)\} \neq \emptyset$,

since $f$ is a neutrosophic continuous function. (142)

3. Let $f(x)$ and $g(x)$ be two neutrosophic mereo-continuous functions at $x = c$, where $f, g: A \to B$. Then,

$$(f+g)(x), (f-g)(x), (f \cdot g)(x), \left(\frac{f}{g}\right)(x) \qquad (143)$$

are all neutrosophic mereo-continuous functions at $x = c$.

### Proofs

$f(x)$ is mereo-continuous at $x = c$ it means that

$$\left\{\lim_{\substack{x\to c\\ x<c}} f(x)\right\} \cap \left\{\lim_{\substack{x\to c\\ x>c}} f(x)\right\} \cap \{f(c)\} \neq \emptyset \qquad (144)$$

therefore

$$\left\{\lim_{\substack{x\to c\\ x<c}} f(x)\right\} = M_1 \cup L_1 \qquad (145)$$

$$\left\{\lim_{\substack{x\to c\\ x>c}} f(x)\right\} = M_1 \cup R_1 \qquad (146)$$

and

$$\{f(c)\} = M_1 \cup V_1 \qquad (147)$$

where all $M_1, L_1, R_1, V_1$ are subsets of $B$, and $M_1 \neq \emptyset$, while $L_1 \cap R_1 \cap V_1 = \emptyset$.

Similarly, $g(x)$ is mereo-continuous at $x = c$ means that

$$\left\{\lim_{\substack{x\to c\\ x<c}} g(x)\right\} \cap \left\{\lim_{\substack{x\to c\\ x>c}} g(x)\right\} \cap \{g(c)\} \neq \emptyset, \qquad (148)$$

therefore

$$\left\{\lim_{\substack{x\to c\\ x<c}} g(x)\right\} = M_2 \cup L_2 \qquad (149)$$

$$\left\{\lim_{\substack{x\to c\\ x>c}} g(x)\right\} = M_2 \cup R_2 \qquad (150)$$

and

$$\{g(c)\} = M_2 \cup V_2 \qquad (151)$$

where all $M_2, L_2, R_2, V_2$ are subsets of $B$, and $M_2 \neq \emptyset$, while $L_2 \cap R_2 \cap V_2 = \emptyset$.

Now,

$f + g: A \to B$

$$(f + g)(x) = f(x) + g(x) \qquad (152)$$

and $(f+g)(x)$ is mereo-continuous at $x=c$ if

$$\left\{\lim_{\substack{x\to c\\ x<c}}(f+g)\right\}\cap\left\{\lim_{\substack{x\to c\\ x>c}}(f+g)(x)\right\}\cap\{(f+g)(c)\}\neq\emptyset \tag{153}$$

or

$$\left\{\lim_{\substack{x\to c\\ x<c}}[f(x)+g(x)]\right\}\cap\left\{\lim_{\substack{x\to c\\ x>c}}[f(x)+g(x)]\right\}\cap \{f(c)+g(c)\}\neq\emptyset \tag{154}$$

or

$$\left(\left\{\lim_{\substack{x\to c\\ x<c}}f(x)\right\}+\left\{\lim_{\substack{x\to c\\ x<c}}g(x)\right\}\right)\cap\left(\left\{\lim_{\substack{x\to c\\ x>c}}f(x)\right\}+\left\{\lim_{\substack{x\to c\\ x>c}}g(x)\right\}\right)\cap(\{f(c)\}+\{g(c)\})\neq\emptyset \tag{155}$$

or

$$(M_1\cup L_1+M_2\cup L_2)\cap(M_1\cup R_1+M_2\cup R_2)\cap (M_1\cup V_1+M_2\cup V_2)\neq\emptyset. \tag{156}$$

But this intersection is non-empty, because:
if $m_1\in M_1\neq\emptyset$ and $m_2\in M_2\neq\emptyset$,
then $m_1\in M_1\cup L_1$, and $m_1\in M_1\cup R_1$, and $m_1\in M_1\cup V_1$ (*)
and $m_2\in M_2\cup L_2$, and $m_2\in M_2\cup R_2$, and $m_2\in M_2\cup V_2$ (**)
whence $m_1+m_2\in M_1\cup L_1+M_2\cup L_2$,
and $m_1+m_2\in M_1\cup R_1+M_2\cup R_2$,
and $m_1+m_2\in M_1\cup V_1+M_2\cup V_2$.

Therefore $(f+g)(x)$ is also mereo-neutrosophic function at $x=c$.

Analogously, one can prove that $f-g$, $f \cdot g$ and $\frac{f}{g}$ are neutrosophic mereo-continuous functions at $x=c$.

From above, one has:

$m_1 - m_2 \in M_1 \cup L_1 - M_2 \cup L_2;$ (157)

$m_1 - m_2 \in M_1 \cup R_1 - M_2 \cup R_2;$ (158)

$m_1 - m_2 \in M_1 \cup V_1 - M_2 \cup V_2.$ (159)

therefore $(f-g)(x)$ is a neutrosophic mereo-continuous function at $x=c$.

Again, from above one has:

$m_1 \cdot m_2 \in (M_1 \cup L_1) \cdot (M_2 \cup L_2);.$ (160)

$m_1 \cdot m_2 \in (M_1 \cup R_1) \cdot (M_2 \cup R_2);.$ (161)

$m_1 \cdot m_2 \in (M_1 \cup V_1) \cdot (M_2 \cup V_2).$ (162)

therefore $(f \cdot g)(x)$ is a neutrosophic mereo-continuous function at $x=c$.

And, from (*) and (**) one has:

$\frac{m_1}{m_2} \in \frac{M_1 \cup L_1}{M_2 \cup L_2};$ (163)

$\frac{m_1}{m_2} \in \frac{M_1 \cup R_1}{M_2 \cup R_2};.$ (164)

$\frac{m_1}{m_2} \in \frac{M_1 \cup V_1}{M_2 \cup V_2}..$ (165)

therefore $\left(\frac{f}{g}\right)(x)$ is a neutrosophic mereo-continuous function at $x=c$.

# III.16. Properties of Neutrosophic Continuity

Similarly to the classical calculus, if $f(x), g(x)$ are neutrosophic continuous functions at $x = c$, and $\alpha \in \mathbb{R}$ is a scalar, then $\alpha \cdot f(x), (f+g)(x), (f-g)(x), (fg)x$ , and $\left(\frac{f}{g}\right)x$ for $g(x) \neq c$ are neutrosophic continuous functions at $x = c$.

The proofs are straightforward as in classical calculus.

Since $f(x)$ and $g(x)$ are neutrosophic continuous functions, one has:

$$\lim_{\substack{x \to c \\ x<c}} f(x) \equiv \lim_{\substack{x \to c \\ x>c}} f(x) \equiv f(c) \tag{166}$$

and
$$\lim_{\substack{x \to c \\ x<c}} g(x) \equiv \lim_{\substack{x \to c \\ x>c}} g(x) \equiv g(c) \tag{167}$$

1. If we multiply the relation (166) by $\alpha$ we get:

$$\alpha \cdot \lim_{\substack{x \to c \\ x<c}} f(x) \equiv \alpha \cdot \lim_{\substack{x \to c \\ x>c}} f(x) \equiv \alpha \cdot f(c) \tag{168}$$

or

$$\lim_{\substack{x \to c \\ x<c}} [\alpha \cdot f(x)] \equiv \lim_{\substack{x \to c \\ x>c}} [\alpha \cdot f(x)] \equiv \alpha \cdot f(c) \tag{169}$$

or $\alpha \cdot f(x)$ is a neutrosophic continuous function at $x = c$.

2. If we add relations (166) and (167) term by term, we get:

$$\lim_{\substack{x \to c \\ x<c}} f(x) + \lim_{\substack{x \to c \\ x<c}} g(x) \equiv \lim_{\substack{x \to c \\ x>c}} f(x) + \lim_{\substack{x \to c \\ x>c}} g(x) \equiv f(c) + g(c) \tag{170}$$

or

$$\lim_{\substack{x\to c\\ x<c}}[f(x)+g(x)] \equiv \lim_{\substack{x\to c\\ x>c}}[f(x)+g(x)] \equiv f(c)+g(c)$$

(171)

or $(f+g)(x)$ is a neutrosophic continuous function at $x = c$.

3. Similarly, if we subtract relations (#) and (##) term by term, we get:

$$\lim_{\substack{x\to c\\ x<c}} f(x) - \lim_{\substack{x\to c\\ x<c}} g(x) \equiv \lim_{\substack{x\to c\\ x>c}} f(x) - \lim_{\substack{x\to c\\ x>c}} g(x) \equiv f(c)-g(c)$$

(172)

or

$$\lim_{\substack{x\to c\\ x<c}}[f(x)-g(x)] \equiv \lim_{\substack{x\to c\\ x>c}}[f(x)-g(x)] \equiv f(c)-g(c)$$

(173)

or $(f-g)(x)$ is a neutrosophic continuous function at $x = c$.

4. If we multiply relations (#) and (##) term by term, we get:

$$\left[\lim_{\substack{x\to c\\ x<c}} f(x)\right]\cdot\left[\lim_{\substack{x\to c\\ x<c}} g(x)\right] \equiv \left[\lim_{\substack{x\to c\\ x>c}} f(x)\right]\cdot\left[\lim_{\substack{x\to c\\ x>c}} g(x)\right] \equiv f(c)\cdot g(c)$$

(174)

or

$$\lim_{\substack{x\to c\\ x<c}}[f(x)\cdot g(x)] \equiv \lim_{\substack{x\to c\\ x>c}}[f(x)\cdot g(x)] \equiv f(c)\cdot g(c)$$

(175)

or $(f\cdot g)(x)$ is a neutrosophic continuous function at $x = c$.

5. If we divide relations (#) and (##) term by term, supposing $g(x) \neq 0$ for all $x$, we get:

$$\frac{\lim_{\substack{x\to c\\ x<c}} f(x)}{\lim_{\substack{x\to c\\ x<c}} g(x)} \equiv \frac{\lim_{\substack{x\to c\\ x>c}} f(x)}{\lim_{\substack{x\to c\\ x>c}} g(x)} \equiv \frac{f(c)}{g(c)} \tag{176}$$

or

$$\lim_{\substack{x\to c\\ x<c}} \left[\frac{f(x)}{g(x)}\right] \equiv \lim_{\substack{x\to c\\ x>c}} \left[\frac{f(x)}{g(x)}\right] \equiv \frac{f(c)}{g(c)} \tag{177}$$

or $\left(\frac{f}{g}\right)(x)$ is a neutrosophic continuous function at $x = c$.

# III.17. The M-δ Definitions of the Neutrosophic Infinite Limits

*The* $M-\delta$ *definitions of the neutrosophic infinite limits* are extensions of the classical infinite limits.

a. $\lim_{x\to c} f(x) = +\infty$ means that $\forall M > 0, \exists \delta = \delta(M) > 0$, such that if $\eta(x,c) < \delta$, then $inf\{f(x)\} > M$.

b. $\lim_{x\to c} f(x) = -\infty$ means that $\forall N < 0, \exists \delta = \delta(N) > 0$, such that if $\eta(x,c) < \delta$, then $sup\{f(x)\} < N$.

# III.18. Examples of Neutrosophic Infinite Limits

1. Let's have the neutrosophic function $f(x) = \frac{[2,\ 5]}{x-1}$.

$$\lim_{\substack{x \to 1 \\ x<1}} \frac{[2,5]}{x-1} = -\infty \tag{178}$$

and

$$\lim_{\substack{x \to 1 \\ x>1}} \frac{[2,5]}{x-1} = +\infty\,. \tag{179}$$

Therefore, $x = 1$ is a vertical asymptote for $f(x)$.

Let's apply the definition for the neutrosophic left limit.

Let $N < 0$. If, for $x < 1$,

$$\eta(x, c) = \eta(x, 1) = |x - 1| < \frac{[2,5]}{|N|} = \delta(N) = \delta, \tag{180}$$

which is equivalent to

$$-\frac{[2,5]}{|N|} < x - 1 < \frac{[2,5]}{|N|} \tag{181}$$

then

$$f(x) = \frac{[2,5]}{x-1} < \frac{[2,5]}{-\frac{[2,5]}{|N|}} = -|N| = N \tag{182}$$

Therefore,

$$\lim_{\substack{x \to 1 \\ x<1}} f(x) = -\infty \tag{183}$$

2. Let $(x) = \frac{4}{(1,3)x^2}$.

$$\lim_{\substack{x \to 0 \\ x<0}} \frac{4}{(1,3)x^2} = +\infty \tag{184}$$

and

$$\lim_{\substack{x \to 0 \\ x>0}} \frac{4}{(1,3)x^2} = +\infty, \tag{185}$$

hence

$$\lim_{x \to 0} \frac{4}{(1,3)x^2} = +\infty. \tag{186}$$

Therefore $x = 0$ is a vertical asymptote for the neutrosophic function $g(x)$.

Let's apply the $M - \delta$ definition to compute the same limit.

Let $M > 0$. If

$$\eta(x, c) = \eta(x, 0) = \eta(x) = |x| < \frac{1}{(\sqrt{1}, \sqrt{3})\sqrt{M}} = \delta(m) = \delta \tag{187}$$

then

$$g(x) = \frac{4}{(1,3)x^2} > \frac{4}{(1,3)\cdot\left[\frac{1}{(\sqrt{1},\sqrt{3})\sqrt{M}}\right]^2} = \frac{4}{(1,3)\cdot\frac{1}{(1,3)M}} =$$

$\frac{4}{\frac{(1,3)/(1,3)}{M}} = 4M/(\frac{1}{3}, 3) =$

because (1,3)/(1,3) = (1/3, 3/1) = (1/3, 3)

$$= (\tfrac{4}{3}\text{M}, 12\text{M}) = \text{M}(\tfrac{4}{3}, 12), \text{ and } \mathit{inf}\{\text{M}(\tfrac{4}{3}, 12)\} = \tfrac{4}{3}M > M. \tag{188}$$

Therefore,

$$\lim_{x \to 0} g(x) = +\infty. \tag{189}$$

2. Let $h(x) = \frac{x^2+7}{x-(\text{either 2 or 3})}$ (190)

be a neutrosophic function [meaning that we are not sure if it is $x-2$ or $x-3$], which is actually equivalent to either the classical function $h_1(x) = \frac{x^2+7}{x-2}$ or to the classical function $h_1(x) = \frac{x^2+7}{x-3}$. (191)

Thus,

$$\lim_{\substack{x \to \text{either 2 or 3} \\ x < \text{either 2 or 3 respectively}}} \frac{x^2+7}{x-(\text{either 2 or 3})} = -\infty$$

(192)

and

$$\lim_{\substack{x \to \text{either 2 or 3} \\ x > \text{either 2 or 3 respectively}}} \frac{x^2+7}{x-(\text{either 2 or 3})} = +\infty$$

(193)

Therefore, either $x = 2$ or $x = 3$ is a vertical asymptote for $h(x)$.

5. Another type of neutrosophic limit:

$$\lim_{x \to 2+2I} \frac{x^2+(1+I)x}{2x+4-6I}$$
$$= \frac{(2+3I)^2+(1+I)(2+3I)}{2(2+3I)+4-6I}$$
$$= \frac{4+12I+9I^2+2+3I+2I+3I^2}{4+6I+4-6I}$$
$$= \frac{6+17I+12I^2}{8} = \frac{6+17I+12I}{8} = \frac{6+29I}{8}$$
$$= \frac{6}{8} + \frac{29}{8}I,$$

where $I$ = indeterminacy with $0 \cdot I = 0$ and $I^2 = I$. (194)

# III.19. Set-Argument Set-Values Function

$f: \mathcal{P}(M) \to \mathcal{P}(N), f(A) = B,$ (195)

where $M$ and $N$ are sets, $A \in \mathcal{P}(M)$ or $A \subseteq M$, and $B \in \mathcal{P}(N)$ or $B \subseteq N$.

This is *a generalization of the interval-argument interval-valued function*.

*Example:*

$f: \mathcal{P}(R) \to \mathcal{P}(R)$ (196)

$f(\{1, 3, 5\}) = \{2, 6\}$ (197)

$f([1, 4]) = [2, 3]$ (198)

$f\big((1, 0)\big) = 5$ (199)

$f([-2,\ 3) \cup \{6\}) = x^2 = [4, 9) \cup \{36\}.$ (200)

$\mathcal{P}(M)$ is the set of all subsets of M, and $\mathcal{P}(N)$ is the set of all subsets of $N$.

The partial-metric $\eta$ and the norm $\mu$ are very well defined on $\mathcal{P}(M)$ and $\mathcal{P}(N)$, and the definitions of neutrosophic limit, neutrosophic continuity, neutrosophic derivative, and neutrosophic integral are extensions from classical calculus definitions by using the partial-metric $\eta$ and/or the norm $\mu$.

# III.20. Neutrosophic Derivative

The general definition of the *neutrosophic derivative* of function $f_N(x)$ is:

$$f_N'(X) = \lim_{\mu(H)\to 0} \frac{<\inf f(X+H)-\inf f(X), \sup f(X+H)-\sup f(X)>}{H}. \quad (201)$$

where *<a, b>* means any kind of open / closed / half open-closed interval.

As particular definitions for the cases when H is an interval one has:

$$f_N'(X) = \lim_{[\inf H,\ \sup H]\to[0,\ 0]} \frac{[\inf f(X+H)-\inf f(X), \sup f(X+H)-\sup f(X)]}{[\inf H,\ \sup H]} \quad (202)$$

is the neutrosophic derivative of $f(X)$.

In a simplified way, one has:

$$f_N'(X) = \lim_{h\to 0} \frac{[\inf f(X+h)-\inf f(X), \sup f(X+h)-\sup f(X)]}{h}. \quad (203)$$

Both definitions above are generalizations of the classical derivative definition, since for crisp functions and crisp variables one has:

$$[\inf H, \sup H] \equiv h \quad (204)$$

and $$\inf f(X+H) \equiv \sup f(x+H) \equiv f(x+h) \quad (205)$$

$$\inf f(X) \equiv \sup f(X) \equiv f(x). \quad (206)$$

Let's see some examples:

1) $f(X) = [2x^3 + 7x, x^5]$. (207)

$$f_N'(X) = \lim_{h\to 0} \frac{[2(x+h)^3 + 7(x+h) - 2x^3 - 7x, (x+h)^5 - x^5]}{h} = \left[\lim_{h\to 0} \frac{2(x+h)^3 + 7(x+h - 2x^3 - 7x}{h}, \lim_{h\to 0} \frac{(x+h)^5 - x^5}{h}\right] = \left[\frac{d}{dx}(2x^3 + 7x), \frac{d}{dx}(x^5)\right] = [6x^2 + 7, 5x^4]. \quad (208)$$

2) Let $g: R \to \mathcal{P}(R)$, by

$$g(x) = \begin{cases} [f_1(x), f_2(x)], \text{if } x \leq 5; \\ [f_3(x), f_4(x)], \text{if } x > 5. \end{cases} \quad (209)$$

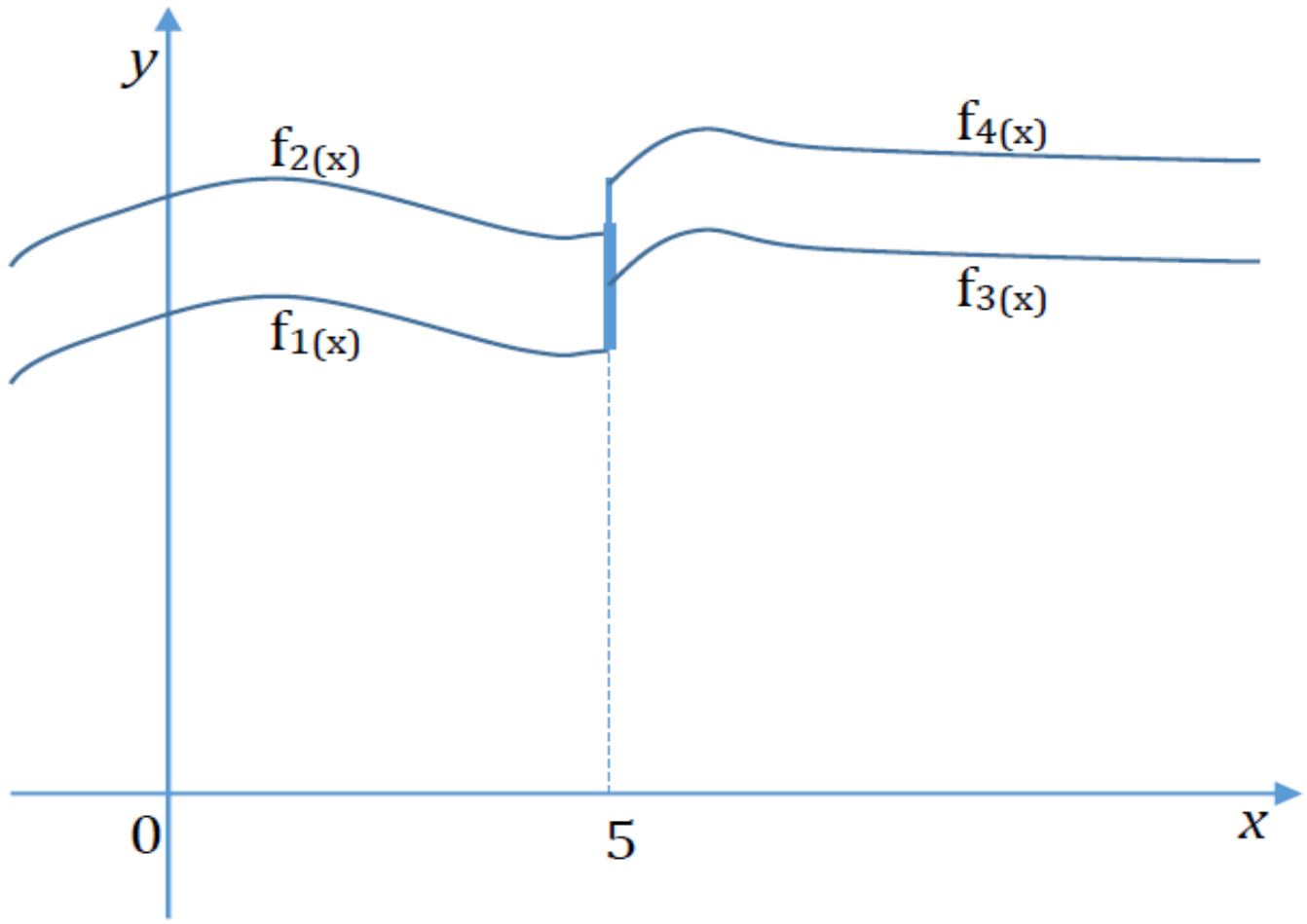

*Graph 23.*

A classical function is differentiable at a given point $x = c$ if: $f$ is continuous at $x = c$, $f$ is smooth at $x = c$, and $f$ does not have a vertical tangent at $x = c$.

$g(x)$ is neutrosophically differentiable on $\mathbb{R} \setminus \{5\}$ if f1, f2, f3, and f4 are differentiable:

$$g'(x) = \begin{cases} [f'_1(x), f'_2(x)], \text{if } x < 5; \\ [f'_3(x), f'_4(x)], \text{if } x > 5. \end{cases} \quad (210)$$

At $x = 5$, the neutrosophic function $g(x)$ is differentiable if:

$$[f'_1(5), f'_2(5)] \equiv [f'_3(5), f'_4(5)], \quad (211)$$

otherwise $g(x)$ has a mereo-derivative at $x = 5$ (as in the above figure) if

$$[f'_1(5), f'_2(5)] \cap [f'_3(5), f'_4(5)] \neq \emptyset, \quad (212)$$

or $g(x)$ is not differentiable at $x = 5$ if

$$[f'_1(5), f'_2(5)] \cap [f'_3(5), f'_4(5)] = \emptyset. \quad (213)$$

3) Another example of neutrosophic derivative.

Let $f : \mathbb{R} \to \mathbb{R} \cup \{I\}$, where $I =$ indeterminacy,

$$f(x) = 3x - x^2 I \quad (214)$$

$$f'(x) = \lim_{h \to 0} \frac{f(x+h) - f(x)}{h}$$

$$= \lim_{h \to 0} \frac{[3(x+h) - (x+h)^2 I] - [3x - x^2 I]}{h}$$

$$= \lim_{h \to 0} \frac{3x + 3h - x^2 I - 2xhI - h^2 I - 3x + x^2 I}{h}$$

$$= \lim_{h \to 0} \frac{h(3 - 2xI - hI)}{h} = 3 - 2xI - 0 \cdot I = 3 - 2xI. \quad (215)$$

Therefore, directly

$$f'(x) = \frac{d}{dx}(3x) - \frac{d}{dx}(x^2 I) = 3 - I\frac{d}{dx}(x^2) = 3 - 2xI. \quad (216)$$

4) An example with refined indeterminacy:

$I_1$ = indeterminacy of first type;

$I_2$ = indeterminacy of second type.

Let $g: \mathbb{R} \to \mathbb{R} \cup \{I_1\} \cup \{I_2\}$, (217)

$g(x) = -x + 2xI_1 + 5x^3I_2$, (218)

Then $g'(x) = \frac{d}{dx}(-x) + \frac{d}{dx}(2xI_1) + \frac{d}{dx}(5x^3I_2) = -1 + 2I_1 + 15x^2I_2$. (219)

# III.21. Neutrosophic Indefinite Integral

We just extend the classical definition of anti-derivative.

The *neutrosophic antiderivative* of neutrosophic function $f(x)$ is the neutrosophic function $F(x)$ such that

$F'(x) = f(x)$.

For example,

1. Let $f: R \to R \cup \{I\}, f(x) = 5x^2 + (3x+1)I$. (220)

Then,

$$F(X) = \int [5x^2 + (3x+1)I]dx$$
$$= \int 5x^2 dx + \int (3x+1)Idx$$
$$= 5 \cdot \frac{x^3}{3} + I \int (3x+1)dx = \frac{5x^3}{3} + \left(\frac{3x^2}{2} + x\right)I + C, \tag{221}$$

where $C$ is an *indeterminate real constant* (i.e. constant of the form *a+bI*, where *a, b* are real numbers, while I = indeterminacy).

2. Refined Indeterminacy.

Let $g: \mathbb{R} \to \mathbb{R} \cup \{I_1\} \cup \{I_2\} \cup \{I_3\}$, (222)

were $I_1, I_2$, and $I_3$ are types of subindeterminacies,

$g(x) = -5 + 2I_1 - x^4 I_2 + 7xI_3$. (223)

Then,
$\int g(x)dx = \int[-5 + 2I_1 - x^4 I_2 + 7xI_3]dx = -5x + 2xI_1 - \frac{x^5}{5}I_2 + \frac{7x^2}{2}I_3 + a + bI$, where $a$ and $b$ are real constants. (224)

# III.22. Neutrosophic Definite Integral

1. Let $h: \mathbb{R} \to \mathcal{P}(\mathbb{R})$ (225)

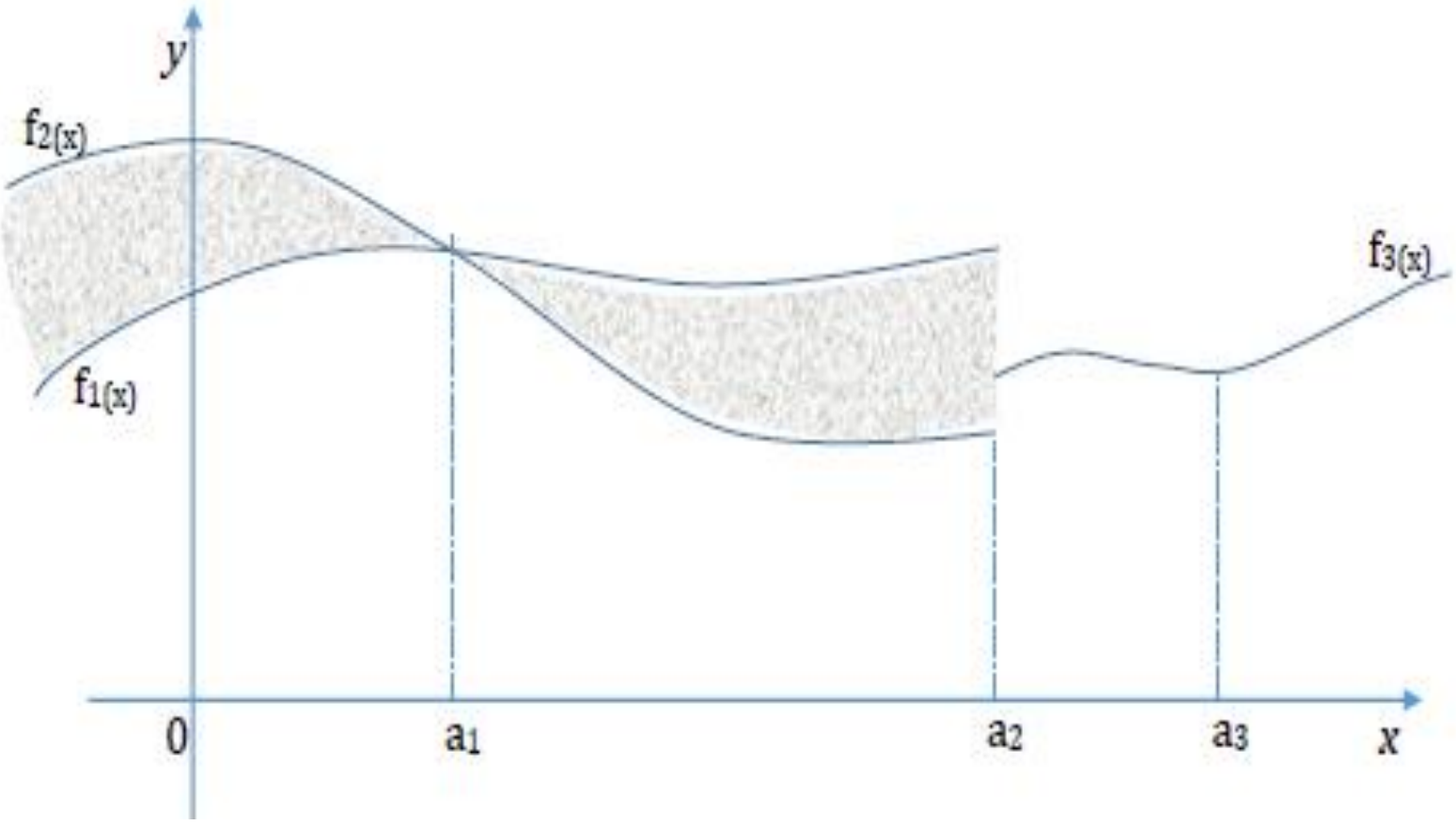

*Graph 24.*

such that

$$h(x) = \begin{cases} [f_1(x), f_2(x)], \text{if } x \leq a_2 \\ f_3(x), \text{if } a > a_2 \end{cases}. \quad (226)$$

$h(x)$ is a thick neutrosophic function for $x \in (-\infty, a_2]$, and a classical function for $x \in (a_2, +\infty)$.

We now compute the *neutrosophic definite integral*:

$\alpha = \int_0^{a_3} h(x)dx = \int_0^{a_1}[f_1(x), f_2(x)]dx + \int_{a_1}^{a_2}[f_2(x), f_1(x)]dx + \int_{a_2}^{a_3} f(x)dx = \left[\int_0^{a_1} f_1(x)dx, \int_0^{a_1} f_2(x)dx\right] + \left[\int_{a_1}^{a_2} f_2(x)dx, \int_{a_1}^{a_2} f_1(x)dx\right] + \int_{a_2}^{a_3} f_3(x)dx = [A, B] + [C, D] + [E, E] = [A + B + E, B + D + E],$ (227)

where, of course,

$A = \int_0^{a_1} f_1(x)dx$, $B = \int_0^{a_1} f_2(x)dx$, $C = \int_{a_1}^{a_2} f_2(x)dx$, $D = \int_{a_1}^{a_2} f_1(x)dx$, and $C = \int_{a_3}^{a_3} f_3(x)dx$.

(228)

Since $h(x)$ is a thick function between *0* and $a_2$, we interpret the result $\alpha$ of our neutrosophic definite integral in general as:

$$\alpha \in [A + B + E, B + D + E], \quad (229)$$

since one may take: $\alpha = A + B + E$ as in classical calculus (i.e. the area are below the lowest curve), or an average:

$$\alpha = \frac{(A + B + E) + (B + D + E)}{2} = \frac{A + D}{2} + B + E \quad (230)$$

(i.e. the area below a curve passing through the middle of the shaded area), or the maximum possible area:

$$\alpha = B + D + E. \quad (231)$$

Depending on the problem to solve, a neutrosophic expert can choose the most appropriate

$$\alpha \in [A + B + E, B + D + E]. \quad (232)$$

# III.23. Simple Definition of Neutrosophic Definite Integral

Let $f_N$ be a neutrosophic function

$$f_N\colon \mathbb{R} \to \mathcal{P}(\mathbb{R}) \tag{233}$$

which is continuous or mereo-continous on the interval $[a, b]$. Then,

$$\Sigma_a^b f_N(x)dx = \lim_{n\to\infty} \Sigma_{i=1}^{n} f_N(C_i)\frac{b-a}{n} \tag{234}$$

where $C_i \in [x_{i-1}, x_i]$, for $i \in \{1, 2, \dots, n\}$, and $a \equiv x_0 < x_1 < x_2 < \cdots < x_{n-1} < x_n \equiv b$ are subdivision of the interval $[a, b]$: exactly as the definition of the classical integral, but $f_N(C_i)$ may be a real set (not necessarily a crisp real number as in classical calculus), or $f_N(C_i)$ may have some indeterminacy.

# III.24. General Definition of Neutrosophic Definite Integral

Let

$f_N: \mathcal{P}(M), \to \mathcal{P}(N),$ (235)

where $M, N$ are given sets, and $\mathcal{P}(M)$ and $\mathcal{P}(N)$ are the power sets of $M$ and $N$ respectively.

$f_N$ is a set-argument set-valued function which, in addition, has some indeterminacy. So, $f_N$ is a neutrosophic set-argument set-valued function.

$f_N$ maps a set in $M$ into a set in $N$. Therefore, $A, B \in \mathcal{P}(M)$. Then:

$$\int_A^B f_N(x)dx = \lim_{n\to\infty} \sum_{i=1}^n f_N(C_i) \cdot \frac{\eta(B,A)}{n}, \quad (236)$$

where

$\inf A \equiv \inf x_0 < \inf x_1 < \cdots < \inf x_{n-1} < \inf x_n \equiv \inf B$

$\sup A \equiv \sup x_0 < \sup x_1 < \cdots < \sup x_{n-1} < \sup x_n \equiv \sup B$

and $(C_i) \in \mathcal{P}(M)$ such that:

$\inf X_{i-1} \leq \inf C_i \leq \inf X_i$

and

$\sup X_{i-1} \leq \sup C_i \leq \sup X_i$, for $i \in \{1, 2, \dots, n\}$.

Therefore, the neutrosophic integral *lower and upper limits* are sets (not necessarily crisp numbers as in classical calculus), $C_i$, for all $i \in \{1, 2, \dots, n\}$, and similarly $f_N(C_i)$ are sets (not crisp numbers as in classical calculus). And, in addition, there may be some indeterminacy as well with respect to their values.

# IV. Conclusion

Neutrosophic Analysis is a generalization of Set Analysis, which in its turn is a generalization of Interval Analysis.

Neutrosophic Precalculus is referred to indeterminate staticity, while Neutrosophic Calculus is the mathematics of indeterminate change.

The Neutrosophic Precalculus and Neutrosophic Calculus can be developed in many ways, depending on the types of indeterminacy one has and on the methods used to deal with such indeterminacy.

We introduce for the first time the notions of *neutrosophic mereo-limit, neutrosophic mereo-continuity* (in a different way from the classical semi-continuity)*, neutrosophic mereo-derivative* and *neutrosophic mereo-integral* (both in different ways from the fractional calculus), besides the classical definitions of limit, continuity, derivative, and integral respectively.

Future research can be done in neutrosophic fractional calculus.

In this book, we present a few examples of indeterminacies and several methods to deal with these specific indeterminacies, but many other indeterminacies there exist in our everyday life, and they have to be studied and resolved using similar of different methods. Therefore, more research should to be done in the field of neutrosophics.

# V. References

## Published Papers and Books

## Other Articles on Neutrosophics

## Presentations to International Conferences or Seminars

[1] F. Smarandache, Foundations of Neutrosophic set and Logic and Their Applications to Information Fusion, Okayama University of Science, Kroumov Laboratory, Department of Intelligence Engineering, Okayama, Japan, 17 December 2013.

[2] Jean Dezert & Florentin Smarandache, Advances and Applications of Dezert-Smarandache Theory (DSmT) for Information Fusion, presented by F. Smarandache, Osaka University, Department of Engineering Science, Inuiguchi Laboratory, Japan, 10 January 2014.

[3] F. Smarandache, Foundations of Neutrosophic Set and Logic and Their Applications to Information Fusion, Osaka University, Inuiguchi Laboratory, Department of Engineering Science, Osaka, Japan, 10 January 2014.

[4] F. Smarandache, Alpha-Discounting Method for Multicriteria Decision Making, Osaka University,

## Ph. D. Dissertations

Neutrosophic Analysis is a generalization of Set Analysis, which in its turn is a generalization of Interval Analysis.

Neutrosophic Precalculus is referred to indeterminate staticity, while Neutrosophic Calculus is the mathematics of indeterminate change.

The Neutrosophic Precalculus and Neutrosophic Calculus can be developed in many ways, depending on the types of indeterminacy one has and on the methods used to deal with such indeterminacy.

In this book, the author presents a few examples of indeterminacies and several methods to deal with these specific indeterminacies, but many other indeterminacies there exist in our everyday life, and they have to be studied and resolved using similar of different methods. Therefore, more research should to be done in the field of neutrosophics.

The author introduces for the first time the notions of neutrosophic mereo-limit, neutrosophic mereo-continuity (in a different way from the classical semi-continuity), neutrosophic mereo-derivative and neutrosophic mereo-integral (both in different ways from the fractional calculus), besides the classical definitions of limit, continuity, derivative, and integral respectively. Future research may be done in the neutrosophic fractional calculus.